\documentclass[11pt,english,a4paper]{smfart}
\usepackage[english]{babel}

\usepackage[utf8]{inputenc}
\usepackage[T1]{fontenc}
\usepackage{kpfonts}
\usepackage{booktabs}
\usepackage{siunitx}
\usepackage{pgfplots}
\usepackage{floatrow}
\usepackage{caption}
\usepackage[nospace]{varioref}
\usepackage[colorlinks=true,linkcolor=blue,citecolor=blue,urlcolor=red, hypertexnames=true]{hyperref}

\usepackage{csquotes}
\usepackage[backend=biber,maxcitenames=5,maxbibnames=99,style=numeric,sorting=nyt,backref=true]{biblatex}
\usepackage{xcolor}
\usepackage{pbsi}
\usepackage{stmaryrd}
\usepackage{dsfont}
\usepackage[mathscr]{euscript}
\usepackage{cancel}
\usepackage{array}
\usepackage{float}
\usepackage{placeins}

\usepackage{enumerate}

\usepackage{amsmath,amssymb,bm,amsfonts}
\usepackage{nicematrix}
\usepackage[all]{xy}
\entrymodifiers={+!!<0pt,\fontdimen22\textfont2>}

\usepackage{mathrsfs}

\newtheorem{theorem}{Theorem}[section]
\newtheorem{lemma}[theorem]{Lemma}
\newtheorem{observation}[theorem]{Observation}
\newtheorem{corollary}[theorem]{Corollary}
\newtheorem{proposition}[theorem]{Proposition}

\theoremstyle{remark}

\newtheorem*{remark*}{\it Remark}

\newcommand{\flba}[2]{
\xymatrix@C15pt{#1\ar@{|->}[r]&#2}}
\newcommand{\flcourte}[2]{
\xymatrix@C12pt{#1\ar[r]&#2}}

{\theoremstyle{definition}}

{\theoremstyle{definition}}

\theoremstyle{definition}

{\theoremstyle{definition}\newtheorem{remark}[theorem]{Remark}}

\def\D{\ensuremath{\mathbb D}}
\def\K{\ensuremath{\mathbb K}}
\def\T{\ensuremath{\mathbb T}}
\def\R{\ensuremath{\mathbb R}}
\def\Z{\ensuremath{\mathbb Z}}

\def\C{\ensuremath{\mathbb C}}
\def\Q{\ensuremath{\mathbb Q}}
\def\N{\ensuremath{\mathbb N}}

\def\bth{\begin{theorem}}
	\def\blm{\begin{lemma}}
		\def\bpr{\begin{proposition}}
			\def\bpf{\begin{proof}}
				\def\epf{\end{proof}}
			\def\epr{\end{proposition}}
		\def\elm{\end{lemma}}
	\def\eth{\end{theorem}}
\def\bco{\begin{corollary}}
	\def\eco{\end{corollary}}
\def\be{\begin{enumerate}}
	\def\ee{\end{enumerate}}

			\newcommand{\vertiii}[1]{{\left\vert\kern-0.25ex\left\vert\kern-0.25ex\left\vert #1 
					\right\vert\kern-0.25ex\right\vert\kern-0.25ex\right\vert}}

			\newcommand{\bx}{{\mathcal B}(X)}

			\newcommand{\lan}{\langle}
			\newcommand{\ran}{\rangle}

			\newcommand{\omg}{\omega}
			\newcommand{\Omg}{\Omega}
			
			\newcommand{\indic}{\mathds{1}}

\numberwithin{equation}{section}

\author[Valentin Gillet]{Valentin Gillet}
\address[Valentin Gillet]{Univ. Lille, CNRS, UMR 8524 - Laboratoire Paul Painlevé, F-59000 Lille, France}
\email{valentin.gillet@univ-lille.fr}

\thanks{This work was supported in part by the project COMOP of the French National Research Agency (grant ANR-24-CE40-0892-01) and by the Labex CEMPI (ANR-11-LABX-0007-01). The author also acknowledges the support of the CDP C2EMPI, as well as the French State under the France-2030 programme, the University of Lille, the Initiative of Excellence of the University of Lille, and the European Metropolis of Lille for their funding and support of the R-CDP-24-004-C2EMPI project.}

\begin{document}

\title[Linear dynamics of random products of weighted shifts]{Linear dynamics of random products of weighted shifts}


\begin{abstract} 
The aim of this article is to study the dynamics of random products of weighted shifts on a separable Fréchet sequence space. That is, given a measure-preserving dynamical system $(\Omg, \mathcal{F}, \mu, \tau) $, a Fréchet sequence space $X$ with a basis $(e_n)_{n \geq 0}$, and a strongly measurable map $T : \Omg \to \bx$ taking values in a finite set of weighted shifts on $X$, we study the dynamics of the sequence $(T(\tau^{n-1} \omg) \dotsm T(\tau \omg) T(\omg))_{n \geq 1}$ for almost every $\omg \in \Omg$. After proving criteria to determine whether this sequence is universal, weakly mixing or mixing for almost every $\omg \in \Omg$, we study some examples on the spaces $X = \ell_p$, $X = c_0$ and $X = H(\C)$ involving two shifts, first in the commuting case and then in the non-commuting one.
\end{abstract}

\maketitle

\par\bigskip
\section{Introduction}\label{Section introduction}

The linear dynamics of random products of operators was studied in \cite{Gill3} for special kinds of operators, such as adjoints of multiplication operators on the Hardy space $H^2(\D)$ or functions of exponential type of the derivation operator on the space of entire functions $H(\C)$. This line of research is motivated by the study of the asymptotic behavior of random products of matrices, initiated by Bellman in \cite{Bell} and by Furstenberg and Kesten in \cite{FK}. Such products arise naturally, for instance, in the study of differential or difference equations with random coefficients.

\smallskip

On the other hand, dynamical properties of linear cocycles have also been investigated in \cite{Bayart-Costakis-Hadjiloucas10}, \cite{Costakis-Hadjiloucas08} and \cite{Martinez-Mendez-Peris26}. 
The results presented in \cite{Gill3} mainly concern random products of commuting operators. In this article, we aim to further investigate the linear dynamics of random products in the non-commutative setting. In particular, we consider random products of weighted shifts, which provide a natural class of examples in this setting. Our approach is to characterize the linear dynamics of these random products through criteria involving the weights of the shifts composing the products. In particular, we show that the universality of these random products can be related to the asymptotic behavior of certain weighted Birkhoff sums.

\smallskip

The random products considered in this article are also related to questions arising in ergodic theory; see, e.g., \cite{Froyland-Lloyd-Quas13} and \cite{Gonzalez-Quas14}.

\medskip

In this paper, we study the linear dynamics of random products of the form
\begin{align} \label{typeproduitsaléaintro}
	T_n(\omg) = T(\tau^{n-1}\omg) \dotsm T(\tau \omg) T(\omg),
\end{align}
where $(A_1, \dotsc, A_k)$ is a fixed partition of measurable subsets of a probability space $\Omg$, and where the operator $T(\omg)$ on $X$ equals a fixed continuous weighted shift operator $B_{w^{(j)}}$ on a separable Fréchet sequence space $X$ when $\omg \in A_j$ for every $1 \leq j \leq k$. 

\smallskip

More precisely, let $(\Omg, \mathcal{F}, \mu, \tau) $ be a measure-preserving dynamical system. Let $(A_1,\dotsc, A_k)$ be a partition of measurable subsets of $\Omg$ with $\mu(A_l) > 0$ for every $1 \leq l \leq k$. We consider weights $w^{(l)} = (w^{(l)}_n)_{n \geq 1}$ for $1 \leq l \leq k$ (with $w^{(l)}_n \ne 0$ for every $1 \leq l \leq k$ and $n \geq 1$), and the corresponding weighted shift operators $B_{w^{(l)}}$ on a separable Fréchet sequence space $X$ given by $B_{w^{(l)}} e_0 = 0$ and $B_{w^{(l)}} e_n = w_n^{(l)} e_{n-1}$ for every $n \geq 1$, where $(e_n)_{n \geq 0}$ is the canonical basis of $X$. We set 
$$T(\omg) = B_{w^{(l)}} \quad \textrm{for every $\omg \in \Omg$ when $\omg \in A_l$, $1 \leq l \leq k$},$$ 
and we study the linear dynamics of the sequence $(T_n(\omg))_{n \geq 1}$ given by (\ref{typeproduitsaléaintro}) for almost every $\omg \in \Omg$. The goal of this article is to find conditions on the weights for the sequence $(T_n(\omg))_{n \geq 1}$ to be universal, weakly mixing, or mixing for almost every $\omg \in \Omg$. To this aim, we consider the random variables $\varepsilon_n, n \geq 1$, defined by 
\begin{align*}
	\varepsilon_n(\omg) = w_n^{(l)} \quad \textrm{if $\omg \in A_l$ for $1 \leq l \leq k$},
\end{align*}
for every $n \geq 1$ and $\omg \in \Omg$.

\subsection{Organization of the paper}
The first part of the article, which corresponds to Section \ref{sectioncriteriaproduitsweightedshifts}, is devoted to proving a characterization of universality, weak mixing, and mixing for the sequence $(T_n(\omg))_{n \geq 1}$ in terms of the weights composing the random products. In Section \ref{sectioncommutingshifts}, we study the dynamics of the sequence $(T_n(\omg))_{n \geq 1}$ generated by two commuting weighted shifts on a separable Fréchet sequence space. Finally, we study in Section \ref{sectiontwoexamplesofnoncommutingshifts} two examples of random products for which the weighted shifts do not commute. The aim of these examples is to illustrate that the dynamics of random products of weighted shifts is related to certain weighted Birkhoff sums and may be delicate in certain situations.

\subsection{Notation and definition} \label{sectionnotetdefinitions}

In this subsection, we gather some important definitions and notation needed for the understanding of this article.

\medskip

\paragraph*{{\textbf{Basic ergodic facts}}}
Let $(\Omg, \mathcal{F}, \mu)$ be a probability space. A transformation $\tau : \Omg \to \Omg$ is measure-preserving if $\mu(\tau^{-1}(B)) = \mu(B)$ for every measurable subset $B$ of $\Omg$. A measure-preserving transformation $\tau $ is said to be ergodic on $(\Omg,\mu)$ if for every measurable subset $A$ of $\Omg$ such that $\tau^{-1}(A) = A$, we have $\mu(A) = 0$ or $\mu(A) = 1$.
 
\smallskip

Two important examples of ergodic transformations are given by irrational rotations and the doubling map on the torus $\T := \R / \Z$. Given an irrational number $\alpha \in (0,1)$, the rotation of parameter $\alpha$ is given by $R_\alpha x = x + \alpha$ in $\T$. The doubling map is the transformation $\tau$ given by $\tau x = 2x$ in $\T$. 

\smallskip
We will also make use of Birkhoff's ergodic theorem in this article, and we recall it here.

\smallskip

\bth \label{thbirkhoffintro}
Let $(\Omg, \mathcal{F}, \mu, \tau)$ be an ergodic measure-preserving dynamical system and $f : \Omg \to \R$ be an integrable function. Then
$$
\frac{1}{n} \displaystyle \sum_{l=0}^{n-1} f(\tau^l \omg) \underset{n\to \infty}{\longrightarrow} \int_\Omg f d\mu \quad \textrm{for almost every $\omg \in \Omg$}.
$$
\eth

\smallskip
Applying Theorem \ref{thbirkhoffintro} to $\indic_{A_l}$, we obtain that $\frac{a_l(n,\omg)}{n} \underset{n\to \infty}{\longrightarrow} \mu(A_l)$ for almost every $\omg \in \Omg$, where $a_l(n,\omg) := \textrm{Card}\{ 0 \leq i \leq n-1 : \tau^i \omg \in A_l \} $.

\medskip

We refer to P. Walters's book \cite{Wal} and to K. Petersen's book \cite{Pet} for background on ergodic theory.

\paragraph*{{\textbf{Linear dynamics}}} We put here the definitions of linear dynamics that will be used in this article.

\smallskip
A sequence of continuous maps $(T_n)_{n \geq 1}$ from a metric space $X$ to itself is said to be universal if there exists an element $x \in X$ such that its orbit under $(T_n)_{n \geq 1}$,
$$
\textrm{orb}(x,(T_n)) := \{ T_n x : n \geq 1 \},
$$
is dense in $X$. Such an $x \in X$ is called a universal vector for $(T_n)_{n \geq 1}$.

A sequence of continuous maps $(T_n)_{n \geq 1}$ from a metric space $X$ to itself is said to be weakly mixing if for every nonempty open sets $U_1, U_2, V_1, V_2$ in $X$, there exists an integer $n \geq 1$ such that $T_n(U_1) \cap V_1 \ne \emptyset$ and $T_n(U_2) \cap V_2 \ne \emptyset$, and it is mixing if for every nonempty open sets $U, V$ in $X$, there exists an integer $N \geq 1$ such that $T_n(U) \cap V \ne \emptyset$ for every $n \geq N$. Mixing obviously implies weak mixing and weak mixing implies universality.

\smallskip

For a single continuous map $T$ on $X$, we say that $T$ is weakly mixing or mixing, if the sequence $(T^n)_{n \geq 1}$ has this property. We say that $T$ is hypercyclic if the sequence $(T^n)_{n \geq 1}$ is universal. 

\smallskip

In the linear setting, we have the following sufficient criterion which implies universality, called the Universality Criterion, and originally proved in \cite[Corollary 1.4]{Godefroy-Shapiro91}.
\smallskip

\bpr[{\cite[Theorem 3.24]{GEP}}] \label{critunivforsequence}
Let $X$ be a separable Fréchet space and $(T_n)_{n \geq 1}$ a sequence of operators on this space. Suppose that there are dense subsets $\mathcal{D}_1 $ and $\mathcal{D}_2$ of $X$, a strictly increasing sequence $(n_k)_{k \geq 1}$ of positive integers, and maps $S_{n_k} : \mathcal{D}_2 \to X$, $k \geq 1$, such that
\begin{enumerate}[(i)]
	\item $T_{n_k} x \underset{k\to \infty}{\longrightarrow} 0\quad \textrm{for every} \; x \in \mathcal{D}_1$.
	\item $S_{n_k} y \underset{k\to \infty}{\longrightarrow} 0\quad \textrm{for every} \; y \in \mathcal{D}_2$.
	\item $T_{n_k} S_{n_k} y \underset{k\to \infty}{\longrightarrow} y \quad \textrm{for every} \; y \in \mathcal{D}_2$.
\end{enumerate}
Then the sequence $(T_n)_{n \geq 1}$ is topologically weakly mixing, and in particular universal. If moreover this criterion is satisfied with respect to the full sequence $(n)$, then the sequence $(T_n)_{n \geq 1}$ is topologically mixing. 
\epr

\smallskip

The maps $S_n$ involved in Proposition \ref{critunivforsequence} are not required to be linear nor continuous. The Universality Criterion will play an important role in our study.  For more details on linear dynamics, we refer to the books \cite{BM} and \cite{GEP}.

\paragraph*{{\textbf{Fréchet sequence spaces}}} By a Fréchet space $X$, we mean a vector space $X$, endowed with a separating increasing sequence $(p_n)_{n \geq 1}$ of seminorms, which is complete for the metric
$$
d(x,y) = \displaystyle \sum_{n \geq 1} 2^{-n} \min(1, p_n(x-y)).
$$
The separating property means that if $p_n(x) = 0$ for every $n \geq 1$, then $x = 0$.
We define an F-norm on $X$ by setting
$$\lVert x \rVert = \displaystyle \sum_{n \geq 1} 2^{-n} \min(1, p_n(x))
$$
for every $x \in X$. One only needs to be careful about the fact that $\lVert . \rVert$ does not satisfy positive homogeneity. However, we have that $\lVert \lambda x \rVert \leq \lVert x \rVert$ for every $x \in X$ and for every $\lambda \in \K$ such that $\lvert \lambda \rvert \leq 1$, and we have that $\lVert \lambda x \rVert \leq (1 + \lvert \lambda \rvert) \lVert x \rVert$ for every $x \in X$ and for every $\lambda \in \K$. We refer to \cite[Section 2.1]{GEP} for more details.

\smallskip
A sequence space $X$ is a linear subspace of $\K^{\Z_+}$ for which the embedding $X \to \K^{\Z_+} $ is continuous, where $\K \in \{\R, \C\}$. In particular, the convergence in $X$ implies the coordinatewise convergence. A Banach/Fréchet sequence space is a Banach/Fréchet space of this kind. For $n \geq 0$, we denote by $e_n$ the sequence given by $e_n(j) = 1$ if $n = j$ and $e_n(j) = 0$ otherwise. For a vector $x \in X$, we write $x_n = \lan e_n^*, x \ran$ for every $n \geq 0$.

\smallskip
Classical examples of Fréchet sequence spaces are of course $X = \K^{\Z_+}$, $ X = \ell_p$ with $1 \leq p < \infty$ or $X = c_0$ with the canonical basis $(e_n)_{n \geq 0}$, as well as $X = H(\C)$, the space of entire functions. In the latter case, an entire function $f(z) = \displaystyle \sum_{n \geq 0} a_n z^n$ is identified with the sequence $(a_n)_{n \geq 0}$, and this sequence space is given by $\{ (a_n)_{n \geq 0} , \displaystyle \lim_{n \to \infty} \lvert a_n \rvert^{1/n} = 0 \} = \{ (a_n)_{n \geq 0} , \displaystyle \sum_{n \geq 0} \lvert a_n \rvert m^n < \infty \quad \textrm{for every}\; m \geq 1 \} $. For every $n \geq 0$, the sequence $e_n$ corresponds to the monomial $z \mapsto z^n$.
\smallskip

Let $X$ be a separable Fréchet sequence space. Suppose that each $e_n$, $n \geq 0$, belongs to $X$ and that $\textrm{span}\{ e_n : n \geq 0\}$ is dense in $X$. If $(w_n)_{n \geq 1}$ is a sequence of \textbf{nonzero} complex numbers, we define the weighted shift operator $B_w$ by
\begin{align}
	B_w e_0 = 0 \quad \textrm{and} \quad B_w e_n = w_n e_{n-1} \quad \textrm{for every } n \geq 1.
\end{align}
When $B_w(X) \subset X$, every weighted shift $B_w$ is continuous on $X$ thanks to closed graph theorem. On $\K^{\Z_+}$, every weighted shift is continuous. On $X = \ell_p$ or $X = c_0$, a weighted shift $B_w$ is continuous if and only if the sequence $w$ is bounded. And finally, on $X = H(\C)$, a weighted shift $B_w$ is continuous if and only if $\displaystyle \sup_{n \geq 1} \lvert w_n \rvert^{1/n} < \infty$. 

\smallskip

The dynamics of a weighted shift $B_w : X \to X$ is well-known in the case where $(e_n)_{n \geq 0}$ is a basis of $X$: $B_w$ is hypercyclic if and only if $B_w$ is weakly mixing, which is also equivalent to the fact that there exists a strictly increasing sequence of positive integers $(n_j)_{j \geq 1}$ such that $$\left(\displaystyle \prod_{l=1}^{n_j} w_l \right)^{-1} e_{n_j} \underset{j\to \infty}{\longrightarrow} 0 \quad \textrm{in $X$.} $$ 
And it is mixing if and only if the same property holds with respect to the full sequence $(n)$. On $X = \ell_p$ or $X = c_0$, the operator $B_w$ is hypercyclic/weakly mixing if and only if $$\displaystyle \sup_{n \geq 1} \prod_{l=1}^n \lvert w_l \rvert = \infty, $$ and it is mixing if and only if $$\displaystyle \lim_{n \to \infty} \prod_{l=1}^n \lvert w_l \rvert = \infty. $$ 
On $X = H(\C)$, the operator $B_w$ is hypercyclic/weakly mixing if and only if $$\displaystyle \sup_{n \geq 1} \left( \prod_{l=1}^n \lvert w_l \rvert \right)^{1/n} = \infty, $$ 
and it is mixing if and only if $$\displaystyle \lim_{n \to \infty} \left(\prod_{l=1}^n \lvert w_l \rvert \right)^{1/n} = \infty. $$ 
Finally, every weighted shift $B_w$ on $\K^{\Z_+}$ is mixing.

\smallskip
For more details and examples on the notion of Fréchet sequence spaces, we refer to \cite[Section 4.1]{GEP}.

\section{Criteria for the universality of random products of weighted shifts} \label{sectioncriteriaproduitsweightedshifts}

In this section, we denote by $X$ a separable Fréchet sequence space, and we denote by $(e_n)_{n \geq 0}$ the sequence of canonical units of $X$, where $\lan e_m^* ,e_n \ran = \delta_{n,m}$ for every $n,m \geq 0$. We suppose that $(e_n)_{n \geq 0}$ is a basis of $X$. 
\smallskip

For a weight $w = (w_n)_{n \geq 1}$, we denote by $B_w$ the weighted shift operator on $X$, given by $B_w e_0 = 0$ and $B_w e_n = w_n e_{n-1}$ for every $n \geq 1$. We consider a measure-preserving dynamical system $(\Omg, \mathcal{F}, \mu, \tau)$, as well as a measurable partition $(A_1, \dotsc, A_k)$ of $\Omg$ with $\mu(A_j) > 0$ for every $1 \leq j \leq k$. We consider weights $w^{(l)} = (w_n^{(l)})_{n \geq 1}$ such that the corresponding weighted shifts $B_{w^{(l)}}$ are continuous on $X$. We give necessary and sufficient conditions on the weights for the sequence $(T_n(\omg))_{n \geq 1}$ to be universal, weakly mixing or mixing for $\omg \in \Omg$, when 
$$
T(\omg) = B_{w^{(l)}} \quad \textrm{if $\omg \in A_l$, for every $\omg \in \Omg$ and $1 \leq l \leq k$}.
$$

For every $n \geq 1$ and $\omg \in \Omg$, we set 
$$
\varepsilon_n(\omg) = w_n^{(l)} \quad \textrm{if $\omg \in A_l$ for $1 \leq l \leq k$}.
$$  
The maps $\varepsilon_n$, $n \geq 1$, define random variables on $\Omg$. Let us remark the following fact.

\smallskip

\begin{observation} \label{observej} For every $n \geq 1$,
	 \begin{align}
	\label{observeqT_n1}	T_n(\omg) e_j &= 0 \quad  \textrm{if \, $0 \leq j < n$} \\
	\label{observeqT_n2}	T_n(\omg) e_j &= \varepsilon_{j-n+1}(\tau^{n-1}\omg) \dotsm \varepsilon_{j-1}(\tau \omg) \varepsilon_j(\omg) \, e_{j-n} \quad \textrm{if \, $j \geq n$}.
	\end{align} 
\end{observation}

\smallskip

With Observation \ref{observej}, we can prove the following result.

\smallskip

\bpr \label{propcondinecessweaklixingFréchet}
Let $X$ be a separable Fréchet sequence space in which $(e_n)_{n \geq 0}$ is a basis and let $(\Omg, \mathcal{F}, \mu, \tau)$ be a measure-preserving dynamical system.
Suppose that the sequence $(T_n(\omg))_{n \geq 1}$ is universal for some $\omg \in \Omg$. Then, there exists a strictly increasing sequence $(n_i)_{i \geq 1}$ of positive integers (depending on $\omg$) such that for every $j \geq 0$,
\begin{align*}
	\left(\displaystyle \prod_{l=1}^{n_i +j} \varepsilon_l(\tau^{n_i +j -l} \omg) \right)^{-1} e_{n_i +j} \underset{i \to \infty}{\longrightarrow} 0 \quad \textrm{in $X$}.
\end{align*} 
\epr

\bpf
Let $\lVert . \rVert$ be an F-norm on $X$ associated to the increasing and separating sequence $(p_n)_{n \geq 1}$ of seminorms.
Let $N \geq 1$ and let $\varepsilon > 0$. We will prove that there exists an integer $n \geq N$ such that for every $0 \leq j \leq N$,
$$\left \lVert \left(\displaystyle \prod_{l=j+1}^{n +j} \varepsilon_l(\tau^{n +j -l} \omg) \right)^{-1} e_{n +j} \right \rVert < \varepsilon.
$$
Since $x_n e_n \underset{n\to \infty}{\longrightarrow} 0$  for every $x \in X$, there exists $\delta > 0$ such that for every $x \in X$,
\begin{equation} \label{eq1prv1}
	\lVert x \rVert < \delta \implies \lVert x_n e_n \rVert < \varepsilon/2 \quad \textrm{for every $n \geq 0$,}
\end{equation}
by the Banach-Steinhaus theorem. Moreover, the convergence in $X$ implies coordinatwise convergence, so there exists $\eta > 0$ such that for every $x \in X$,
\begin{equation} \label{eq2prv1}
	\lVert x \rVert <\eta \implies \lvert \lan e_l^*, x \ran \rvert \leq 1/2 \quad \textrm{for every $0 \leq l \leq N$}.
\end{equation}
Let $y \in X$ be a universal vector for $(T_n(\omg))_{n \geq 1}$. Let us notice that $ry$ remains a universal vector for $(T_n(\omg))_{n \geq 1}$, for every $r > 0$. Indeed, let $v \in X, m \geq 1$ and $\alpha > 0$. There exists an integer $n \geq 1$ such that $p_m(T_n(\omg) y - \frac{v}{r}) < \alpha/r$, and thus $p_m(T_n(\omg)(ry) - v) < \alpha$, since $p_m$ is a seminorm.

We now show that there exists a vector $x \in X$ with $\lVert x \rVert < \delta$ and an integer $n \geq N$ such that $\lVert T_n(\omg) x - \displaystyle \sum_{0 \leq l \leq N} e_l \rVert < \eta$. Let $N' \geq 1$ be an integer such that $p_{N'}(y) > 0$ and $\displaystyle \sum_{n > N'} 2^{-n} < \delta/2$. We choose $r > 0$ such that $p_{N'}(ry) = r p_{N'}(y) < \delta/2$. Using the fact that sequence of seminorms is increasing, we obtain that
\begin{align*}
	\lVert r y \rVert &\leq \displaystyle \sum_{n=1}^{N'} 2^{-n} p_n(ry) + \displaystyle \sum_{n > N'} 2^{-n} \\
	&< \delta/2 + p_{N'}(ry) \\
	&< \delta.
\end{align*}
We now set $x := r y$, which satisfies $\lVert x \rVert < \delta$. Since $X$ has no isolated point,  the set $\{ T_n(\omg) x : n \geq N \}$ remains dense in $X$ and there exists an integer $n \geq N$ such that $\lVert T_n(\omg) x - \displaystyle \sum_{0 \leq l \leq N} e_l \rVert < \eta$, as desired. 

By (\ref{eq2prv1}), we obtain that
$$
\lvert \lan e_j^*, T_n(\omg) x - \displaystyle \sum_{0 \leq l \leq N} e_l \ran \rvert \leq 1/2 \quad \textrm{for every $0 \leq j \leq N$,}
$$
that is, 
$$
\lvert (\varepsilon_{n+j}(\omg) \varepsilon_{n+j-1}(\tau \omg) \dotsm \varepsilon_{j+1}(\tau^{n-1}\omg)) x_{n +j} -1  \rvert \leq 1/2
$$
for every $0 \leq j \leq N$,
by (\ref{observeqT_n2}).
In particular, we obtain that
$$
\left \lvert \frac{(\varepsilon_{n+j}(\omg) \varepsilon_{n+j-1}(\tau \omg) \dotsm  \varepsilon_{j+1}(\tau^{n-1}\omg)) x_{n+j} -1 }{(\varepsilon_{n+j}(\omg) \varepsilon_{n+j-1}(\tau \omg) \dotsm \varepsilon_{j+1}(\tau^{n-1} \omg)) x_{n+j}} \right \rvert \leq 1
$$
for every $0 \leq j \leq N$,
and it follows from (\ref{eq1prv1}) that
\begin{align*}
	&\left \lVert \left(\displaystyle \prod_{l=j+1}^{n+j} \varepsilon_l(\tau^{n+j -l} \omg) \right)^{-1} e_{n+j} \right \rVert \\
	\leq &\left \lVert \left ( \left (\displaystyle \prod_{l=j+1}^{n+j} \varepsilon_l(\tau^{n+j -l} \omg) \right)^{-1} x_{n+j}^{-1} -1 \right) x_{n+j}  e_{n+j} \right \rVert + \lVert x_{n+j} e_{n+j} \rVert \\
	<& \varepsilon
\end{align*}
for every $0 \leq j \leq N$.
\smallskip
In particular, there exists a strictly increasing sequence of positive integers $(n_i)_{i \geq 1}$ (depending on $\omg \in \Omg$) such that for every $j \geq 0$,
\begin{align*}
	\left(\displaystyle \prod_{l=j+1}^{n_i +j} \varepsilon_l(\tau^{n_i +j -l} \omg) \right)^{-1} e_{n_i +j} \underset{i \to \infty}{\longrightarrow} 0 \quad \textrm{in $X$}.
\end{align*} It follows that for every $j \geq 0$,
\begin{align*}
	&\left(\displaystyle \prod_{l=1}^{n_i +j} \varepsilon_l(\tau^{n_i +j -l} \omg) \right)^{-1} e_{n_i +j}  \\
	=& \left(\displaystyle \prod_{l=1}^{j} \varepsilon_l(\tau^{n_i +j -l} \omg) \right)^{-1} \left(\displaystyle \prod_{l=j+1}^{n_i +j} \varepsilon_l(\tau^{n_i +j -l} \omg) \right)^{-1} e_{n_i +j}
\end{align*} 
converges to $0$ in $X$ as $i \to \infty$, since the products
$$
\left(\displaystyle \prod_{l=1}^{j} \lvert \varepsilon_l(\tau^{n_i +j -l} \omg) \rvert \right)^{-1}
$$
can be bounded from above independently of $i \geq 1$. This concludes the proof of Proposition \ref{propcondinecessweaklixingFréchet}.
\epf

\smallskip

We now prove the converse of Proposition \ref{propcondinecessweaklixingFréchet}.

\smallskip

\bpr \label{propcondisuffisanteweaklmixingFréchet}
Let $X$ be a Fréchet sequence space in which $(e_n)_{n \geq 0}$ is a basis and let $(\Omg, \mathcal{F}, \mu, \tau)$ be a measure-preserving dynamical system.

Let $\omg \in \Omg$. Suppose that there exists a strictly increasing sequence of positive integers $(n_i)_{i \geq 1}$ (depending on $\omg$) such that for every $j \geq 0$,
\begin{align*}
	\left(\displaystyle \prod_{l=1}^{n_i +j} \varepsilon_l(\tau^{n_i +j -l} \omg) \right)^{-1} e_{n_i +j} \underset{i \to \infty}{\longrightarrow} 0 \quad \textrm{in $X$}.
\end{align*} 
Then the sequence $(T_n(\omg))_{n \geq 1}$ is weakly mixing.
\epr

\smallskip

\bpf
We apply the Universality Criterion (Proposition \ref{critunivforsequence}). Let $X_0 = \textrm{span}\{ e_n : n \geq 0\}$, which is dense in $X$. We define a right inverse $S_n(\omg)$ of $T_n(\omg)$ by setting 
$$S_n(\omg) e_j = (\varepsilon_{n+j}(\omg) \varepsilon_{n+j-1}(\tau \omg) \dotsm \varepsilon_{j+2}(\tau^{n-2}\omg) \varepsilon_{j+1}(\tau^{n-1}\omg))^{-1} e_{n+j}$$ for every $n \geq 1$ and $j \geq 0$. Let us notice that $T_n(\omg) e_j \underset{n\to \infty}{\longrightarrow} 0 $ in $X$ for every $j \geq 0$. Moreover, since 
\begin{align*}
	S_{n_i}(\omg) e_j &= \left(\displaystyle \prod_{l=j+1}^{n_i +j} \varepsilon_l(\tau^{n_i +j -l} \omg) \right)^{-1} e_{n_i +j} \\
	&= \left(\displaystyle \prod_{l=1}^{j} \varepsilon_l(\tau^{n_i +j -l} \omg) \right) 	\left(\displaystyle \prod_{l=1}^{n_i +j} \varepsilon_l(\tau^{n_i +j -l} \omg) \right)^{-1} e_{n_i +j}
\end{align*}
we have that for every $j \geq 0$, $S_{n_i}(\omg) e_j$ converges to $0$ in $X$ as $i \to \infty$. This concludes the proof of Proposition \ref{propcondisuffisanteweaklmixingFréchet}.  	
\epf

\smallskip

As a consequence of Propositions \ref{propcondinecessweaklixingFréchet} and \ref{propcondisuffisanteweaklmixingFréchet}, we obtain the following result.

\smallskip

\bth \label{theoremcaractweakmixingFréchetBanach}
Let $X$ be a Fréchet sequence space in which $(e_n)_{n \geq 0}$ is a basis and let $(\Omg, \mathcal{F}, \mu,\tau)$ be a measure-preserving dynamical system. Let $\omg \in \Omg$. The following assertions are equivalent.
\begin{enumerate}[(1)]
	\item The sequence $(T_n(\omg))_{n \geq 1}$ is universal;
	\item The sequence $(T_n(\omg))_{n \geq 1}$ is weakly mixing;
	\item There exists a strictly increasing sequence of positive integers $(n_i)_{i \geq 1}$ (depending on $\omg$) such that for every $j \geq 0$,
	\begin{align*}
		\left(\displaystyle \prod_{l=1}^{n_i +j} \varepsilon_l(\tau^{n_i +j -l} \omg) \right)^{-1} e_{n_i +j} \underset{i\to \infty}{\longrightarrow} 0 \quad \textrm{in $X$}.
	\end{align*} 
\end{enumerate}
\eth

\smallskip

The same proof as in Propositions \ref{propcondinecessweaklixingFréchet} and \ref{propcondisuffisanteweaklmixingFréchet} leads to the following result regarding mixing.

\smallskip

\bth \label{thcaractmelangefortprodtransfoinversible}
Let $X$ be a Fréchet sequence space and let $(\Omg, \mathcal{F}, \mu,\tau)$ be a measure-preserving dynamical system. Let $\omg \in \Omg$. The following assertions are equivalent.
\begin{enumerate}[(1)]
	\item The sequence $(T_n(\omg))_{n \geq 1}$ is mixing;
	\item The following condition holds: 
	\begin{align*}
		\left(\displaystyle \prod_{l=1}^{n} \varepsilon_l(\tau^{n -l} \omg) \right)^{-1} e_{n} \underset{n \to \infty}{\longrightarrow} 0 \quad \textrm{in $X$}.
	\end{align*} 
\end{enumerate}
\eth

We end this section with the following remark.

\smallskip

\begin{remark}
	\begin{enumerate}\itemsep0.5em
		\item In the case where $X = \ell_p$ or $X = c_0$ with $1 \leq p < \infty$, the sequence $(T_n(\omg))_{n \geq 1}$ is universal/weakly mixing if and only if there exists a strictly increasing sequence of positive integers $(n_i)_{i \geq 1}$ (depending on $\omg \in \Omg$) such that for every $j \geq 0$,
		\begin{align*}
			\displaystyle \lim_{i \to \infty} \prod_{l=1}^{n_i +j} \lvert \varepsilon_l(\tau^{n_i +j-l}\omg) \rvert = \infty.
		\end{align*}
		
		Moreover, the sequence $(T_n(\omg))_{n \geq 1}$ is mixing if and only if
		\begin{align*}
			\displaystyle \lim_{n \to \infty} \prod_{l=1}^{n} \lvert \varepsilon_l(\tau^{n-l}\omg) \rvert = \infty.
		\end{align*}
		\item In the case where $X = H(\C)$, a sequence $(a_n e_n)_{n \geq 0}$ converges to 0 if and only if $\displaystyle \lim_{n \to \infty} \lvert a_n \rvert^{1/n} = 0$. In particular, the sequence $(T_n(\omg))_{n \geq 1}$ is universal/weakly mixing if and only if there exists a strictly increasing sequence of positive integers $(n_i)_{i \geq 1}$ (depending on $\omg \in \Omg$) such that for every $j \geq 0$,
		\begin{align*}
			\displaystyle \lim_{i \to \infty} \left( \prod_{l=1}^{n_i +j} \lvert \varepsilon_l(\tau^{n_i +j-l}\omg) \rvert \right)^{1/(n_i+j)} = \infty.
		\end{align*}
		
		Moreover, the sequence $(T_n(\omg))_{n \geq 1}$ is mixing if and only if
		\begin{align*}
			\displaystyle \lim_{n \to \infty} \left( \prod_{l=1}^{n} \lvert \varepsilon_l(\tau^{n-l} \omg) \rvert \right)^{1/n} = \infty  
		\end{align*}
		\item In the case where $X = \mathcal{\K}^{\N}$, every weighted shift is continuous on $X$ and the sequence $(T_n(\omg))_{n \geq 1}$ is always mixing for every $\omg \in \Omg$.
		\item In the case where $X = \ell_p$ or $X = c_0$ with $1 \leq p < \infty$, if $(\Omg, \mathcal{F}, \mu,\tau)$ is ergodic and if 
		\begin{align} \label{conditionpoidsnonunivnorminfinie}
			\displaystyle \prod_{l=1}^{k} \lVert w^{(l)} \rVert_\infty^{\mu(A_l)} < 1,
		\end{align}
		then the sequence $(T_n(\omg))_{n \geq 1}$ is not universal for almost every $\omg \in \Omg$. Indeed, suppose that (\ref{conditionpoidsnonunivnorminfinie}) holds. We have that
		\begin{align*}
			\displaystyle \prod_{l=1}^{n} \lvert \varepsilon_l(\tau^{n-l}\omg) \rvert \leq \displaystyle \prod_{l=1}^{k} \lVert w^{(l)} \rVert_\infty^{a_l(n,\omg)}
		\end{align*}
		where
		\begin{align*}
			a_l(n,\omg) = \textrm{Card}\{0 \leq i \leq n-1 : \tau^i \omg \in A_l \}.
		\end{align*}
		By Birkhoff's ergodic theorem, there exists a measurable subset $E$ of $\Omg$ with $\mu(E) = 1$ such that for every $\omg \in E$ and for every $1 \leq l \leq k$,
		\begin{align*}
			\frac{a_l(n,\omg)}{n} \underset{n\to \infty}{\longrightarrow} \mu(A_l).
		\end{align*}
		Since
		\begin{align*}
			\displaystyle \prod_{l=1}^{k} \lVert w^{(l)} \rVert_\infty^{a_l(n,\omg)} = \exp(n \log(\displaystyle \prod_{l=1}^{k} \lVert w^{(l)} \rVert_\infty^{\mu(A_l)}) + o(n))
		\end{align*}
		as $n \to \infty$, the sequence 
		\begin{align*}
			\left(\displaystyle\prod_{l=1}^{n} \lvert \varepsilon_l(\tau^{n-l}\omg) \rvert \right)_{n \geq 1}
		\end{align*}
		converges to $0$ for almost every $\omg \in \Omg$. 
	\end{enumerate}
\end{remark}

\section{The case of two commuting weighted shifts}\label{sectioncommutingshifts}

In this subsection, we focus on the case where $X$ is either $\ell_p, c_0$ or $H(\C)$. We consider the situation of two continuous weighted shifts $B_w$ and $B_v$ on $X$ such that $B_w B_v = B_v B_w$, and 
\[
T(\omega) =
\begin{cases}
	B_w & \text{if } \omega \in A_1,\\
	B_v & \text{if } \omega \in A_2,
\end{cases}
\]
where $(\Omg, \mathcal{F}, \mu, \tau)$ is a measure-preserving dynamical system and  $(A_1,A_2)$ is a measurable partition of $\Omega$ such that $\min\{\mu(A_1), \mu(A_2)\} > 0$. 

\medskip

Let us notice that the weighted shifts $B_w$ and $B_v$ commute if and only if $w_{j+1} v_j = w_j v_{j+1}$ for every $j \geq 1$, that is, if and only if 
\[
w_j = c v_j \quad \text{for every } j \ge 1, \text{ where } c := \frac{w_1}{v_1}.
\]
We will suppose in what follows that the sequences $v$ and $w$ are distinct, that is, $c \ne 1$. In this context, we have that
\begin{align} \label{equationlienvetwcasdeuxshifts}
	\displaystyle \prod_{l=1}^{n} \lvert \varepsilon_l(\tau^{n-l} \omg) \rvert = \lvert c \rvert^{a_1(n,\omg)} \displaystyle \prod_{l=1}^{n} \lvert v_l \rvert = \lvert c \rvert^{-a_2(n,\omg)} \displaystyle \prod_{l=1}^{n} \lvert w_l \rvert
\end{align}
and that
\begin{align} \label{equationlienhcvetwcasdeuxshifts}
	\displaystyle \prod_{l=1}^{n} \lvert w_l \rvert = \lvert c \rvert^{n} \displaystyle \prod_{l=1}^{n} \lvert v_l \rvert
\end{align}
for every $n \geq 1$,
where we recall that $a_l(n,\omg) = \textrm{Card}\{ 0 \leq i \leq n-1 : \tau^i \omg \in A_l \} $.

\smallskip

Before exploring the case where $X = \ell_p$ or $X = c_0$, we do the following observation.

\smallskip
\begin{observation} \label{univBvetsoussuite}
	Let $B_v$ be a weighted shift on $X = \ell_p$ or $X = c_0$, where $v$ is a bounded sequence of nonzero complex numbers. Suppose that there exists a strictly increasing sequence of positive integers $(n_i)_{i \geq 1}$ such that $\displaystyle \prod_{l=1}^{n_i} \lvert v_l \rvert \underset{i\to \infty}{\longrightarrow} \infty$. Then there exists a strictly increasing sequence of positive integers $(m_i)_{i \geq 1}$ such that $\displaystyle \prod_{l=1}^{m_i +j} \lvert v_l \rvert \underset{i\to \infty}{\longrightarrow} \infty$ for every $j \geq 0$.
\end{observation}

\smallskip

\bpf
Let $C > 1$ be a constant such that $\lvert v_l \rvert \leq C$ for every $l \geq 1$.
For every $i \geq 1$, there exists $N_i \geq i+1$ with $N_i \in (n_l)_{l \geq 1}$ such that $\displaystyle \prod_{l=1}^{N_i} \lvert v_l \rvert > C^{2i} $. This implies that for every $0 \leq j \leq i$, we have $\displaystyle \prod_{l=1}^{N_i -j} \lvert v_l \rvert > C^{i}$. We now set $m_i = N_i -i$ for every $i \geq 1$. Then we obtain that $\displaystyle \prod_{l=1}^{m_i + j} \lvert v_l \rvert > C^i$ for every $i \geq 1$ and $0 \leq j \leq i$. Passing to a strictly increasing subsequence of $(m_i)_{i \geq 1}$ if necessary, this concludes the proof of Observation \ref{univBvetsoussuite}.
\epf

\smallskip
With Observation \ref{univBvetsoussuite}, we obtain the following result.

\smallskip

\bpr \label{propcasdeuxshiftscommutentc>1}
Suppose that $X = \ell_p$ or $X = c_0$. Let $(\Omg, \mathcal{F}, \mu, \tau)$ be a measure-preserving dynamical system. Let $(A_1, A_2)$ be a non-trivial measurable partition of $\Omg$ and let $v,w$ be two bounded sequences of nonzero complex numbers. 
Suppose that
\[
T(\omega) =
\begin{cases}
	B_w & \text{if } \omega \in A_1,\\
	B_v & \text{if } \omega \in A_2,
\end{cases}
\] 
with $B_w B_v = B_v B_w$. 

\smallskip

Suppose that $\lvert w_1 \rvert \geq \lvert v_1 \rvert$. Then
\begin{equation*}
	B_v \text{ is hypercyclic} 
	\;\Longrightarrow\; 
	\bigl(T_n(\omega)\bigr)_{n \geq 1} \text{ is universal for some/every } \omega \in \Omega
	\;\Longrightarrow\; 
	B_w \text{ is hypercyclic}
\end{equation*}
and
\begin{equation*}
	B_v \text{ is mixing} 
	\;\Longrightarrow\; 
	\bigl(T_n(\omega)\bigr)_{n \geq 1} \text{ is mixing for some/every } \omega \in \Omega
	\;\Longrightarrow\; 
	B_w \text{ is mixing.}
\end{equation*}
\epr

\smallskip

\bpf
Suppose first that $B_v$ is hypercyclic. Since $\lvert c \rvert \geq 1$, we have that 
\begin{align*}
	\displaystyle \prod_{l=1}^{n +j} \lvert \varepsilon_l(\tau^{n +j-l}\omg) \rvert \geq \prod_{l=1}^{n +j} \lvert v_l \rvert
\end{align*}
for every $n \geq 1$ and $j \geq 0$.
Thus, Observation \ref{univBvetsoussuite} implies that $(T_n(\omg))_{n \geq 1}$ is universal for every $\omg \in \Omg$.

\smallskip

Suppose now that $(T_n(\omg))_{n \geq 1}$ is universal for some $\omg \in \Omg$. Since
\begin{align*}
	\displaystyle \prod_{l=1}^{n} \lvert w_l \rvert = \lvert c \rvert^{n} \displaystyle \prod_{l=1}^{n} \lvert v_l \rvert \geq \displaystyle \prod_{l=1}^{n} \lvert \varepsilon_l(\tau^{n-l} \omg) \rvert,
\end{align*}
we obtain that $\displaystyle \sup_{n \geq 1} \prod_{l=1}^{n} \lvert w_l \rvert = \infty$, that is, $B_w$ is hypercyclic. 

The same arguments show the implications for mixing. This concludes the proof of Proposition \ref{propcasdeuxshiftscommutentc>1}.
\epf

\smallskip

\begin{remark}
In the case where $\lvert w_1 \rvert = \lvert v_1 \rvert$ in Proposition \ref{propcasdeuxshiftscommutentc>1}, we notice that the converse is also true, that is, 
\begin{align*}
	(T_n(\omg))_{n \geq 1 } \text{ is universal for some/every } \omega \in \Omega \iff B_v \text{ is hypercyclic} \iff B_w \text{ is hypercyclic} 
\end{align*}
and
\begin{align*}
	(T_n(\omg))_{n \geq 1 } \text{ is mixing for some/every } \omega \in \Omega \iff B_v \text{ is mixing} \iff B_w \text{ is mixing}. 
\end{align*}
\end{remark}

\medskip

We now move to the case where $\lvert w_1 \rvert \ne \lvert v_1 \rvert$ and show that in this case the converse of Proposition \ref{propcasdeuxshiftscommutentc>1} does not hold in general.

\smallskip

\bpr \label{Propexempledeuxshiftscommutentc>1}
Suppose that $X = \ell_p$ or $X = c_0$. Let $(\Omg, \mathcal{F}, \mu, \tau)$ be an ergodic measure-preserving dynamical system. Let $(A_1, A_2)$ be a non-trivial measurable partition of $\Omg$. 
Suppose that
\[
T(\omega) =
\begin{cases}
	B_w & \text{if } \omega \in A_1,\\
	B_v & \text{if } \omega \in A_2,
\end{cases}
\] 
with $B_w B_v = B_v B_w$. 

Suppose that $\lvert w_1 \rvert > \lvert v_1 \rvert$.
\begin{enumerate}
	\item If $\displaystyle \inf_{n \geq 1} \prod_{l=1}^{n} \lvert v_l \rvert > 0 $, then $(T_n(\omg))_{n \geq 1}$ is mixing for almost every $\omg \in \Omg$.
	\item If $n^{-1} \displaystyle \sum_{l=1}^{n} \log \lvert v_l \rvert \underset{n \to \infty}{\longrightarrow} \beta' $ with $ - \log \lvert c \rvert < \beta' < -\mu(A_1) \log \lvert c \rvert $, then $B_w$ is mixing and $(T_n(\omg))_{n \geq 1}$ is not universal for almost every $\omg \in \Omg$.
\end{enumerate}
\epr

\smallskip

\bpf
The first point of Proposition \ref{Propexempledeuxshiftscommutentc>1} just follows from (\ref{equationlienvetwcasdeuxshifts}) and from the fact that $\lvert c \rvert^{a_1(n,\omg)}  \underset{n \to \infty}{\longrightarrow} \infty $ for almost every $\omg \in \Omg$. 

For the other points, we look at $\log \left(\displaystyle \prod_{l=1}^{n} \lvert \varepsilon_l(\tau^{n-l} \omg) \rvert \right)  = a_1(n,\omg) \log\lvert c \rvert +  \displaystyle \sum_{l=1}^{n} \log \lvert v_l \rvert$.

Suppose that $n^{-1} \displaystyle \sum_{l=1}^{n} \log \lvert v_l \rvert \underset{n \to \infty}{\longrightarrow} \beta' $ with $ - \log \lvert c \rvert < \beta' < -\mu(A_1) \log \lvert c \rvert $. Then, for almost every $\omg \in \Omg$,
\begin{align*}
	a_1(n,\omg) \log \lvert c \rvert + \displaystyle \sum_{l=1}^{n} \log \lvert v_l \rvert \sim n(\mu(A_1) \log \lvert c \rvert + \beta')
\end{align*} 
and
\begin{align*}
	 \log \left(\displaystyle \prod_{l=1}^{n} \lvert w_l \rvert \right) = n \log \lvert c \rvert + \displaystyle \sum_{l=1}^{n} \log \lvert v_l \rvert \sim n( \log \lvert c \rvert + \beta')
\end{align*} 
as $n \to \infty$. Since $-\log \lvert c \rvert < \beta' < -\mu(A_1) \log \lvert c \rvert$, we see that $B_w$ is mixing and that $(T_n(\omg))_{n \geq 1}$ is not universal for almost every $\omg \in \Omg$.
\epf

\smallskip

In general, the situation where $\lvert w_1 \rvert > \lvert v_1 \rvert$
and $\displaystyle \prod_{l=1}^{n} \lvert v_l \rvert
\underset{n \to \infty}{\longrightarrow} 0$
is delicate and it is not always easy to decide whether
$(T_n(\omega))_{n \geq 1}$ is universal for almost every $\omega \in \Omega$. We illustrate it now. It is based on the properties of Birkhoff sums associated to an ergodic measure-preserving transformation and to a centered function. More precisely, given an ergodic measure-preserving dynamical system $(\Omg, \mathcal{F}, \mu, \tau)$ and a centered function $f : \Omg \to \R$, it is known (see, for instance, \cite{Ha}) that the sequence $(\mathbb{S}_n^{\tau}f(\omg))_{n \geq 1}$ of Birkhoff sums $ \mathbb{S}_n^{\tau}f(\omg) := \displaystyle \sum_{l=0}^{n-1} f(\tau^l \omg) $ changes sign infinitely often for almost every $\omg \in \Omg$, and that the sequence $(\mathbb{S}_n^{\tau}f(\omg))_{n \geq 1}$ is bounded for almost every $\omg \in \Omg$ if and only if the coboundary equation $f = h - h \circ \tau$ has a solution $h \in L^\infty(\Omg)$ (see, for instance, \cite[Section 4.1, Theorem 19]{Kachu}).

\smallskip

\begin{proposition} \label{excasdeuxshiftscommutentavecsommesbirkhoff}
	Suppose that $X = \ell_p$ or $X = c_0$. Let $(\Omg, \mathcal{F}, \mu, \tau)$ be an ergodic measure-preserving dynamical system. Let $(A_1, A_2)$ be a non-trivial measurable partition of $\Omg$.
	Suppose that
	\[
	T(\omega) =
	\begin{cases}
		B_w & \text{if } \omega \in A_1,\\
		B_v & \text{if } \omega \in A_2,
	\end{cases}
	\] 
	with $B_w B_v = B_v B_w$.  
	
	\smallskip
	
	Suppose that $\lvert w_1 \rvert > \lvert v_1 \rvert$ and that
	\begin{align*}
		\displaystyle \prod_{l=1}^{n} \lvert v_l \rvert
		\sim \lvert c\rvert^{- n \mu(A_1)}
	\end{align*}
	as $n \to \infty$, where we recall that $c= \frac{w_1}{v_1}$. In this case, $ B_w$ is mixing. Let $f := \indic_{A_1} - \mu(A_1)$ and let $\mathbb{S}_n^{\tau}f(\omg) := \displaystyle \sum_{l=0}^{n-1} f(\tau^l \omg)$ be the $n$-th Birkhoff sum associated to $\tau$ and $f$.
	\begin{enumerate}
		\item If $\displaystyle \limsup_{n \to \infty} \mathbb{S}_n^{\tau}f(\omg) = \infty $ for almost every $\omg \in \Omg$, then $(T_n(\omg))_{n \geq 1}$ is universal for almost every $\omg \in \Omg$;
		\item If $\displaystyle \limsup_{n \to \infty} \mathbb{S}_n^{\tau}f(\omg) < \infty $ for almost every $\omg \in \Omg$, then $(T_n(\omg))_{n \geq 1}$ is not universal for almost every $\omg \in \Omg$.
	\end{enumerate}
\end{proposition}

\smallskip

\bpf
In this case, the equation (\ref{equationlienhcvetwcasdeuxshifts}) implies that 
$$\displaystyle \prod_{l=1}^{n} \lvert w_l \rvert \sim \lvert c \rvert^{n \mu(A_2)}$$ as $n \to \infty$, and in particular $B_w$ is mixing.

Now, for every $j \geq 0$, we have that
\begin{align*}
	\lvert c \rvert^{a_1(n +j,\omg)} \displaystyle \prod_{l=1}^{n+j} \lvert v_l \rvert \sim \lvert c \rvert^{a_1(n,\omg) - n \mu(A_1) } \lvert c \rvert^{-j \mu(A_1) + b_1(n,j,\omg)} \sim \lvert c \rvert^{\mathbb{S}_n^{\tau}f(\omg)} \lvert c \rvert^{-j \mu(A_1) + b_1(n,j,\omg)}
\end{align*} 
as $n \to \infty$, where we set $b_1(n,j,\omg) := \textrm{Card}\{ n \leq i \leq n+j-1 : \tau^i \omg \in A_1\}$. Since $\lvert c \rvert > 1$ and since $\lvert c \rvert^{b_1(n,j,\omg)} \leq \lvert c \rvert^{j}$, the assumptions on $\displaystyle \limsup_{n \to \infty} \mathbb{S}_n^{\tau}f(\omg)$ implies the conclusion of Proposition \ref{excasdeuxshiftscommutentavecsommesbirkhoff}.
\epf

\smallskip
Examples of Birkhoff sums $\mathbb{S}_n^{\tau}f$ satisfying $\displaystyle \limsup_{n \to \infty} \mathbb{S}_n^{\tau}f(\omg) = \infty$ for almost every $\omg \in \Omg$ is given by Birkhoff sums having a subsequence satisfying a Central Limit Theorem (\cite[Proposition 2.5]{Gill3}). Given an ergodic measure-preserving dynamical system $(\Omg, \mathcal{F}, \mu, \tau)$ and a centered function $f : \Omg \to \R$, we say that the sequence $(\mathbb{S}_n^{\tau}f)_{n \geq 1}$ satisfies a Central Limit Theorem (CLT) if there exists a sequence of positive real numbers $(a_n)_{n \geq 1}$ with $a_n\underset{n\to \infty}{\longrightarrow} \infty$ such that the sequence $(\frac{\mathbb{S}_n^{\tau}f}{a_n})_{n \geq 1}$ converges in distribution to $\mathcal{N}(0,1)$, whose density on $\R$ is $x \mapsto \frac{1}{\sqrt{2 \pi}} \exp(-x^2/2) $. An investigation of Birkhoff sums satisfying a CLT for the doubling map and irrational rotations was done in \cite{Gill3} for some special step functions. This investigation was based on the following results. For the doubling map, it is due to Kac.

\smallskip

\bth[\cite{K}] \label{TCLKacdoublangle}
Let $f$ be a centered real function on $\T$ such that the Fourier coefficients of $f$ satisfy
$$
\exists \beta > 1/2,\; \forall n \ne 0, \; \lvert c_n(f) \rvert \leq C/\lvert n \rvert^{\beta}.
$$
Moreover, suppose that the equation $f = g - g \circ \tau$ has no solution $g \in L^2(\T)$, where $\tau$ is the doubling map on $\T$. Then there exists $\sigma^2 > 0$ such that the sequence $(\frac{1}{\sigma \sqrt{n}} \displaystyle \sum_{l=0}^{n-1} f \circ \tau^l )_{n \geq 1}$ converges in distribution to $\mathcal{N}(0,1)$.
\eth

\smallskip

Applying Theorem \ref{TCLKacdoublangle} to the function $f := \indic_{A_1} - m(A_1)$ when $A_1 = [0,b)$ with $b \in (0,1)$ in the context of $\Omg = \T$ and of the doubling map, where $m$ is the normalized Lebesgue measure on $\T$, we obtain for example the following result.

\smallskip 

\medskip

\bco
Suppose that $X = \ell_p$ or $X = c_0$. Let $b \in (0,1)$, $A_1 = [0,b)$, $A_2 = [b,1)$ and let $\Omg = \T$ be equipped with the normalized Lebesgue measure $m$ on $\T$ and with the doubling map $\tau$ on $\T$.
Suppose that
\[
T(\omega) =
\begin{cases}
	B_w & \text{if } \omega \in A_1,\\
	B_v & \text{if } \omega \in A_2,
\end{cases}
\] 
with $B_w B_v = B_v B_w$. 

\smallskip

Suppose that $\lvert w_1 \rvert > \lvert v_1 \rvert$ and that
\begin{align*}
	\displaystyle \prod_{l=1}^{n} \lvert v_l \rvert
	\sim \lvert c\rvert^{- n \mu(A_1)}
\end{align*}
as $n \to \infty$. Then $(T_n(\omg))_{n \geq 1}$ is universal for almost every $\omg \in \T$.
\eco

\smallskip

\bpf
The proof can easily be adapted from the proof given in \cite[Lemma 2.20]{Gill3}.
\epf

\smallskip

For rotations, the problem is more difficult. We will illustrate it with the case where the irrational parameter has bounded partial quotients. For an irrational number $\alpha \in (0,1)$, we denote by $(a_j)_{j \geq 1}$ the sequence of partial quotients of its continued fraction. For a centered function $f : \T \to \R$ of bounded variation, its Fourier coefficients $c_r(f)$ can be written as $c_r(f) = \frac{\gamma_r(f)}{r} $ for every $r \ne 0$, with $\displaystyle \sup_{r \ne 0} \lvert \gamma_r(f) \rvert < \infty$. The following CLT holds for irrational rotations whose parameters have bounded partial quotients.

\smallskip

\bth [{\cite[Theorem 3.6]{CB}}]
Let $f : \T \to \R$ be a centered function of bounded variation and $\alpha$ an irrational number with bounded partial quotients, which means that $\displaystyle \sup_{n \geq 1} a_n < \infty$. Suppose that the following property holds:
$$
\exists \eta > 0,\, \theta > 0,\, M \geq 1, \quad  \operatorname{Card}\{0 \leq j \leq N : a_{j+1} \lvert \gamma_{q_j}(f) \rvert \geq \eta \} \geq \theta N, \; \forall N \geq M,
$$
where $(q_j)_{j \geq 1}$ is the sequence of denominators of the convergents of $\alpha$. Then there exists a subset $W$ of positive integers with density $1$ on which $\lVert \mathbb{S}_n^{R_\alpha}f \rVert_2 \underset{n\to \infty}{\longrightarrow} \infty $ such that the sequence $ \left(\frac{\mathbb{S}_n^{R_\alpha}f}{\lVert \mathbb{S}_n^{R_\alpha}f \rVert_2} \right)_{n \in W}$ converges in distribution to $\mathcal{N}(0,1)$.  
\eth

\medskip

\bco
Suppose that $X = \ell_p$ or $X = c_0$. Let $\alpha$ be an irrational number in $(0,1)$ and let $\Omg = \T$ equipped with the normalized Lebesgue measure $m$ on $\T$ and with the rotation $R_\alpha$. Let $A_1 = [0,b)$ and $A_2 = [b,1)$. Suppose that $\alpha$ has bounded partial quotients.

Suppose that
\[
T(\omega) =
\begin{cases}
	B_w & \text{if } \omega \in A_1,\\
	B_v & \text{if } \omega \in A_2,
\end{cases}
\] 
with $B_w B_v = B_v B_w$. 

\smallskip

Suppose that $\lvert w_1 \rvert > \lvert v_1 \rvert$ and that
\begin{align*}
	\displaystyle \prod_{l=1}^{n} \lvert v_l \rvert
	\sim \lvert c\rvert^{- n \mu(A_1)}
\end{align*}
as $n \to \infty$. If $b$ is a rational number, then $(T_n(\omg))_{n \geq 1}$ is universal for almost every $\omg \in \T$.
\eco

\smallskip

\bpf
The proof can easily be adapted from \cite[Corollary 2.12]{Gill3}.
\epf

\medskip

The case where $X = H(\C)$ is more simple than the case where $X = \ell_p$ or $X = c_0$, thanks to Birkhoff's ergodic theorem.

\smallskip

\bpr \label{excasdeuxshiftscommutentfoncentieres}
Suppose that $X = H(\C)$.  Let $(\Omg, \mathcal{F}, \mu, \tau)$ be an ergodic measure-preserving dynamical system. Let $(A_1, A_2)$ be a non-trivial measurable partition of $\Omg$. 
Suppose that
\[
T(\omega) =
\begin{cases}
	B_w & \text{if } \omega \in A_1,\\
	B_v & \text{if } \omega \in A_2,
\end{cases}
\] 
with $B_w B_v = B_v B_w$.

\smallskip

The following assertions are equivalent regarding universality/weak mixing.
\begin{enumerate}
	\item $(T_n(\omg))_{n \geq 1}$ is universal for almost every $\omg \in \Omg$;
	\item $(T_n(\omg))_{n \geq 1}$ is weakly mixing for almost every $\omg \in \Omg$;
	\item $B_v$ is hypercyclic;
	\item $B_w$ is hypercyclic.
\end{enumerate}

Moreover, regarding mixing, the following assertions are equivalent.
\begin{align*}
	(T_n(\omg))_{n \geq 1} \; \textrm{is mixing for a.e $\omg \in \Omg$} \iff B_v \; \textrm{is mixing} \iff B_w \; \textrm{is mixing}.
\end{align*}
\epr

\smallskip

\bpf
Let us notice that in this case, we have that
$$
\displaystyle \left(\prod_{l=1}^n \lvert w_l \rvert \right)^{1/n} = \lvert c \rvert \displaystyle \left (\prod_{l=1}^n \lvert v_l \rvert \right)^{1/n}
$$
and that
$$
\displaystyle \left (\prod_{l=1}^n \lvert \varepsilon_l(\tau^{n-l} \omg) \rvert \right)^{1/n} = \lvert c \rvert^{\frac{a_1(n,\omg)}{n}}  \left(\displaystyle \prod_{l=1}^n \lvert v_l \rvert \right)^{1/n}
$$
for every $n \geq 1$, by (\ref{equationlienvetwcasdeuxshifts}) and (\ref{equationlienhcvetwcasdeuxshifts}).

By Birkhoff's ergodic theorem, the sequence $(\frac{a_1(n,\omg)}{n})_{n \geq 1}$ converges to $\mu(A_1)$ for almost every $\omg \in \Omg$. Thus if $(T_n(\omg))_{n \geq 1}$ is universal/weakly mixing for almost every $\omg \in \Omg$, we obtain that $B_v$ and $B_w$ are hypercyclic. Conversely, if $B_v$ is hypercyclic, then there exists a strictly increasing sequence of positive integers $(n_i)_{i \geq 1}$ such that $\displaystyle \left (\prod_{l=1}^{n_i} \lvert v_l \rvert \right)^{1/{n_i}} \underset{i\to \infty}{\longrightarrow} \infty$. As in Observation \ref{univBvetsoussuite}, this implies that there exists a strictly increasing sequence of positive integers $(m_i)_{i \geq 1}$ such that $\displaystyle \left (\prod_{l=1}^{m_i +j} \lvert v_l \rvert \right)^{1/{(m_i +j)}} \underset{i\to \infty}{\longrightarrow} \infty$ for every $j \geq 0$. This implies that $(T_n(\omg))_{n \geq 1}$ is universal for almost every $\omg \in \Omg$.
\epf

\section{Study of two examples of non-commuting weighted shifts} \label{sectiontwoexamplesofnoncommutingshifts}

In this subsection, we study two examples of random products of weighted shifts on $X = \ell_p$, $X = c_0$ or $X = H(\C)$ that do not commute. These examples rely on weighted ergodic theorems. More precisely, let $X$ be either one of these sequence spaces, let $B_w$ and $B_v$ be two continuous weighted shifts on $X$. Let $(\Omg, \mathcal{F}, \mu, \tau)$ be an ergodic measure-preserving dynamical system. Let $(A_1, A_2)$ be a non-trivial measurable partition of $\Omg$. 
Suppose that
\[
T(\omega) =
\begin{cases}
	B_w & \text{if } \omega \in A_1,\\
	B_v & \text{if } \omega \in A_2.
\end{cases}
\] 
Then for every $n \geq 1$,
\begin{align*}
		\log \displaystyle \left (\prod_{l=1}^{n} \lvert \varepsilon_l(\tau^{n-l} \omg) \rvert \right) = \displaystyle \sum_{l=0}^{n-1} \log(\lvert w_{n-l} \rvert) \indic_{A_1}(\tau^l \omg) + \displaystyle \sum_{l=0}^{n-1} \log(\lvert v_{n-l} \rvert) \indic_{A_2}(\tau^l \omg).
\end{align*}
Normalizing these sums by the numbers $\displaystyle \sum_{l=1}^{n} \log \lvert w_l \rvert$ and $\displaystyle \sum_{l=1}^{n} \log \lvert v_l \rvert$ respectively, we obtain the Nörlund means of the sequence $(f(\tau^l \omg))_{n \geq 1}$ by the weights $w$ and $v$ (for an exposition on Nörlund summation, we refer to Hardy’s book \cite[Section 4.1]{Hardybook}). A study of the convergence of these sums is needed to study the linear dynamics of the sequence $(T_n(\omg))_{n \geq 1}$ for almost every $\omg \in \Omg$ in the case of weighted shifts. While classical Birkhoff sums have been extensively investigated, the corresponding literature on Nörlund-type Birkhoff sums remains rather limited. Our two examples aim to illustrate it.

\subsection{An example on the space $X = \ell_p$ or $X = c_0$}

We consider in this subsection the case where $X = \ell_p$ or $X = c_0$. We consider an ergodic measure-preserving dynamical system $(\Omg, \mathcal{F}, \mu, \tau)$, as well as a non-trivial measurable partition $(A_1,A_2)$ of $\Omg$.

Let 
\[
T(\omega) =
\begin{cases}
	B_w & \text{if } \omega \in A_1,\\
	B_v & \text{if } \omega \in A_2,
\end{cases}
\] 
where the weights are given by 
\begin{align*}
	w_l = 1 + \frac{1}{l+1} \quad \textrm{and} \quad v_l =  1 - \frac{1}{l+1} \quad \textrm{for every $l \geq 1$}. 
\end{align*}

Let us compute $\displaystyle \prod_{l=1}^{n} \lvert \varepsilon_l(\tau^{n-l} \omg) \rvert $ in this case. One has
\begin{align*}
	\displaystyle \prod_{l=1}^{n} \lvert \varepsilon_l(\tau^{n-l} \omg) \rvert &= \displaystyle \prod_{l=1}^{n} \left( \frac{l+2}{l+1} \right)^{\indic_{A_1}(\tau^{n-l} \omg)} \left (\frac{l}{l+1} \right)^{\indic_{A_2}(\tau^{n-l}\omg)} \\
	&= \frac{1}{(n+1)!} \displaystyle \prod_{l=1}^{n} l \left( 1+\frac{2}{l} \right)^{\indic_{A_1}(\tau^{n-l} \omg)} \\
	&= \frac{1}{n+1} \displaystyle \prod_{l=1}^{n} \left (1+ \frac{2}{l} \right)^{\indic_{A_1}(\tau^{n-l}\omg)}.
\end{align*}

Let us set $V_n(\omg) := \log \left(\displaystyle \prod_{l=1}^{n} \lvert \varepsilon_l(\tau^{n-l} \omg) \rvert \right) = -\log(n+1) + \displaystyle \sum_{l=0}^{n-1} \log(1+ \frac{2}{n-l}) \indic_{A_1}(\tau^l \omg)$ for every $n \geq 1$ and $\omg \in \Omg$. The following inequalities hold:
\begin{align} \label{ineqinfVn}
	-\log(n+1)
	&+ 2 \sum_{l=0}^{n-1} \frac{1}{\,n-l\,}\,\indic_{A_1}(\tau^l\omega)
	- 2 \sum_{l=0}^{n-1} \frac{1}{(n-l)^2}\,\indic_{A_1}(\tau^l\omega)
	\le V_n(\omega)
\end{align}
and
\begin{align} \label{ineqsupVn}
	V_n(\omega)
	&\le -\log(n+1)
	+ 2 \sum_{l=0}^{n-1} \frac{1}{\,n-l\,}\,\indic_{A_1}(\tau^l\omega).
\end{align}
Since the sequence $\displaystyle \left (\sum_{l=0}^{n-1} \frac{1}{(n-l)^2}\,\indic_{A_1}(\tau^l\omega) \right)_{n \geq 1}$ is bounded for every $\omg \in \Omg$, we thus have to study the weighted sums $\displaystyle \sum_{l=0}^{n-1} \frac{1}{\,n-l\,}\,\indic_{A_1}(\tau^l\omega)$ for almost every $\omg \in \Omg$.

\medskip

A possible way to treat the universality/weak mixing condition is to use the following results. The first one is a weighted version of Birkhoff's ergodic theorem in $L^2(\Omg)$.

\smallskip

\blm \label{lemmecvgceL2sommesnorlund}
Let $(\Omg, \mathcal{F}, \mu, \tau)$ be an ergodic measure-preserving dynamical system and let $f : \Omg \to \R$ be in $L^2(\Omg)$. Let $(p_n)_{n \geq 1}$ be a sequence of positive real numbers. Suppose that
\begin{align}
	&\frac{1}{P_n} (\lvert p_n - p_{n+1} \rvert + \dotsc + \lvert p_1 - p_{2} \rvert) \underset{n \to \infty}{\longrightarrow} 0, \\
	&\frac{p_{n+1}}{P_n} \underset{n \to \infty}{\longrightarrow} 0 \quad \textrm{and} \\
	&P_n \underset{n \to \infty}{\longrightarrow} \infty,
\end{align}
where $P_n = \displaystyle \sum_{j=1}^n p_j$.
Then \( \displaystyle \frac{1}{P_n} \sum_{l=0}^{n-1} p_{n-l} f(\tau^l \omg) \) converges to $\int_\Omg f d\mu$ in $L^2(\Omg)$.
\elm

\smallskip

Classical examples of sequences $(p_n)_{n \geq 1}$ satisfying Lemma \ref{lemmecvgceL2sommesnorlund} are given by
\begin{enumerate}
	\item $p_n \leq p_{n+1}$ for every $n \geq 1$, $p_{n+1}/P_n \underset{n \to \infty}{\longrightarrow} 0 $ and $P_n \underset{n \to \infty}{\longrightarrow} \infty $;
	\item $p_{n+1} \leq p_n$ for every $n \geq 1$, $p_{n+1}/P_n \underset{n \to \infty}{\longrightarrow} 0 $ and $P_n \underset{n \to \infty}{\longrightarrow} \infty $.
\end{enumerate} 

\smallskip

The proof of Lemma \ref{lemmecvgceL2sommesnorlund} is based on the following observation.

\smallskip

\begin{observation} \label{observdensité}
	If $\tau$ is ergodic on $(\Omg, \mathcal{F}, \mu)$, the set $\{ g - g \circ \tau : g \in L^2(\Omg) \}$ is dense in the set of functions $f$ in $L^2(\Omg)$ such that $\int_\Omg f d\mu = 0$. 
\end{observation}

\smallskip

\bpf[Proof of Observation \ref{observdensité}]
This is well-known. Since the Koopman operator $U_\tau : f \mapsto f \circ \tau$ is an isometry on $L^2(\Omg)$, the set $(U_\tau-I)(L^2(\Omg))$ is dense in $(\textrm{Ker}(U_\tau - I))^\perp$, which is easily seen to be equal to the set $\{ f \in L^2(\Omg) : \int_\Omg f d\mu = 0 \}$, by ergodicity of $\tau$.
\epf

\smallskip

\bpf[Proof of Lemma \ref{lemmecvgceL2sommesnorlund}]
In order to prove Lemma \ref{lemmecvgceL2sommesnorlund}, we can suppose that $\int_\Omg f d\mu = 0$. 

Let $g : \Omg \to \R$ be in $L^2(\Omg)$. Then
\begin{align*}
	\frac{1}{P_n} \displaystyle \sum_{l=0}^{n-1} p_{n-l} (g(\tau^l \omg) - g(\tau^{l+1} \omg)) &= \frac{1}{P_n} \left( \displaystyle \sum_{l=0}^{n-1} (p_{n-l} - p_{n-l + 1}) g(\tau^l \omg) + p_{n+1} g(\omg) - p_1 g(\tau^n \omg) \right)
\end{align*}
and thus
\begin{align*}
	&\left \lVert \frac{1}{P_n} \displaystyle \sum_{l=0}^{n-1} p_{n-l} (g \circ\tau^l - g \circ \tau^{l+1}) \right \rVert_{2} \\
	\leq & \frac{1}{P_n} (p_{n+1} + p_1 + \lvert p_n - p_{n+1} \rvert + \dotsc + \lvert p_2 - p_{1} \rvert ) \lVert g \rVert_{2}
\end{align*}
which converges to 0 as $n \to \infty$.

\smallskip

Now if $f : \Omg \to \R$ is in $L^2(\Omg)$ such that $\int_\Omg f d\mu = 0$ and if $\varepsilon > 0$, then there exists a function $f_\varepsilon$ in $\{ g - g \circ \tau : g \in L^2(\Omg) \}$ such that $\lVert f - f_\varepsilon \rVert_2 < \varepsilon/2$.
Thus 
\begin{align*}
	&\left \lVert \frac{1}{P_n} \displaystyle \sum_{l=0}^{n-1} p_{n-l} f \circ\tau^l \right \rVert_{2} \\
	\leq & \left \lVert \frac{1}{P_n} \displaystyle \sum_{l=0}^{n-1} p_{n-l} (f-f_\varepsilon)\circ\tau^l \right \rVert_{2} + \left \lVert \frac{1}{P_n} \displaystyle \sum_{l=0}^{n-1} p_{n-l} f_\varepsilon \circ\tau^l \right \rVert_{2} \\
	<& \varepsilon
\end{align*}
when $n$ is large enough, since \( \displaystyle \frac{1}{P_n} \sum_{l=0}^{n-1} p_{n-l} f_\varepsilon \circ\tau^l \) converges to 0 in $L^2(\Omg)$ as $n \to \infty$.
This concludes the proof of Lemma \ref{lemmecvgceL2sommesnorlund}.
\epf

\smallskip

The second tool we need is the diagonal extraction principle, given by the following result.

\smallskip

\blm \label{extracdiagonalsuitesL2}
Let $(u_n)_{n \geq 1}$ be a sequence in $L^2(\Omg)$ such that $u_n  \underset{n\to \infty}{\longrightarrow} a$ in $L^2(\Omg)$, with $a \in \R$. Then there exist a strictly increasing map $\varphi : \N \to \N$ and a measurable set $E$ with $\mu(E) = 1$ such that for every $\omg \in E$, $u_{\varphi(n) +j}(\omg)  \underset{n\to \infty}{\longrightarrow} a $ for every $j \geq 0$.
\elm

\smallskip

\bpf
Since $u_n  \underset{n\to \infty}{\longrightarrow} a$ in $L^2(\Omg)$, there exist a strictly increasing map $\varphi_0 : \N \to \N$ and a measurable set $E_0$ with $\mu(E_0) = 1$ such that for every $\omg \in E_0$, $u_{\varphi_0(n)}(\omg) \underset{n\to \infty}{\longrightarrow} a $. Since $u_{\varphi_0(n) +1} \underset{n\to \infty}{\longrightarrow} a$ in $L^2(\Omg)$, there exist a strictly increasing map $\varphi_1 : \N \to \N$ and a measurable set $E_1$ with $\mu(E_1) = 1$ such that for every $\omg \in E_1$, $u_{(\varphi_0 \circ \varphi_1)(n) +1}(\omg) \underset{n\to \infty}{\longrightarrow} a $. Proceeding by induction, we then show that there exist a sequence $(\varphi_n)_{n \geq 0}$ of strictly increasing maps from $\N$ to $\N$, as well as a sequence of measurable sets $(E_n)_{n \geq 0}$ satisfying $\mu(E_n) = 1$ for every $n \geq 0$, such that for every $j \geq 0$ and for every $\omg \in E_j$, we have that $u_{(\varphi_0 \circ \dotsc \circ \varphi_j)(n) + j}(\omg) \underset{n\to \infty}{\longrightarrow} a $. We set $E = \displaystyle \bigcap_{j \geq 0} E_j$ which satisfies $\mu(E) = 1$, and we define $\varphi : \N \to \N$ by $\varphi(n) = (\varphi_0 \circ \dotsc \circ \varphi_n)(n)$ for every $n \geq 1$. For every $j \geq 0$ and for every $\omg \in E$, we have that $u_{\varphi(n) +j}(\omg)  \underset{n\to \infty}{\longrightarrow} a $ and this concludes the proof of Lemma \ref{extracdiagonalsuitesL2}.  
\epf

\smallskip

Using Lemmas \ref{lemmecvgceL2sommesnorlund} and \ref{extracdiagonalsuitesL2}, we can give a condition on $\mu(A_1)$ and $\mu(A_2)$ in order to have universality/weak mixing of the sequence $(T_n(\omg))_{n \geq 1}$ for almost every $\omg \in \Omg$.

\smallskip

\bpr \label{propmélangeexcasdeuxshiftsmesureA_1>1/2et<1/2}
Let $X = \ell_p$ or $X = c_0$.
Suppose that $(\Omg, \mathcal{F}, \mu,\tau)$ is an ergodic measure-preserving dynamical system and that $(A_1,A_2)$ is a partition of two measurable subsets of $\Omega$ with $\min \{\mu(A_1), \mu(A_2) \} > 0$.

Let 
\[
T(\omega) =
\begin{cases}
	B_w & \text{if } \omega \in A_1,\\
	B_v & \text{if } \omega \in A_2,
\end{cases}
\] 
where the weights are given by 
\begin{align*}
	w_l = 1 + \frac{1}{l+1} \quad \textrm{and} \quad v_l =  1 - \frac{1}{l+1} \quad \textrm{for every $l \geq 1$}. 
\end{align*}
\begin{enumerate}
	\item If $\mu(A_1) > \mu(A_2)$, then $(T_n(\omg))_{n \geq 1}$ is universal/weakly mixing for almost every $\omg \in \Omg$;
	\item If $\mu(A_1) < \mu(A_2)$, then $(T_n(\omg))_{n \geq 1}$ is not mixing for almost every $\omg \in \Omg$.
\end{enumerate}
\epr

\smallskip

\bpf
In our case, the weights are given by $p_n = 1/n$ for every $n \geq 1$ and obviously satisfy the conditions of Lemma \ref{lemmecvgceL2sommesnorlund}. In particular, \(\displaystyle \frac{1}{P_n} \sum_{l=0}^{n-1} \frac{1}{n-l} \indic_{A_1} \circ\tau^l \) converges to $\mu(A_1)$ in $L^2(\Omg)$. Let $\varphi : \N \to \N$ be a strictly increasing map and $E$ be a measurable set with $\mu(E) = 1$ such that \(\displaystyle \frac{1}{P_{\varphi(n) +j}} \sum_{l=0}^{\varphi(n) +j-1} \frac{1}{\varphi(n) +j-l} \indic_{A_1} \circ\tau^l(\omg) \) converges to $\mu(A_1)$ for every $\omg \in E$ and for every $j \geq 0$. Since $P_n \sim \log(n+1)$ as $n \to \infty$, we obtain (using (\ref{ineqinfVn})) that 
$$V_{\varphi(n) +j}(\omg) := \log \left(\displaystyle \prod_{l=1}^{\varphi(n) +j} \lvert \varepsilon_l(\tau^{\varphi(n) +j-l} \omg) \rvert \right) \underset{n \to \infty}{\longrightarrow} \infty
$$ 
for every $j \geq 0$ and for every $\omg \in E$, if $\mu(A_1) > 1/2$, which implies that $(T_n(\omg))_{n \geq 1}$ is weakly mixing for almost every $\omg \in \Omg$ when $\mu(A_1) > 1/2$. 

\smallskip

In the case $\mu(A_1) < 1/2$, we obtain (using (\ref{ineqsupVn})) that 
$$V_{\varphi(n)}(\omg) = \log \left(\displaystyle \prod_{l=1}^{\varphi(n)} \lvert \varepsilon_l(\tau^{\varphi(n)-l} \omg) \rvert \right) \underset{n \to \infty}{\longrightarrow} -\infty
$$
for every $\omg \in E$, which implies that $(T_n(\omg))_{n \geq 1}$ is not mixing for almost every $\omg \in \Omg$ when $\mu(A_1) < 1/2$. 
\epf

However, we will see that it is possible to find a non-trivial measurable partition $(A_1, A_2)$ of $\Omg$ with $\mu(A_1) > 1/2$ such that the sequence $(T_n(\omg))_{n \geq 1}$ is weakly mixing but not mixing for almost every $\omg \in \Omg$.

\smallskip

\bpr \label{proppartitionharmonicsumsnotconverge}
Let $X = \ell_p$ or $X = c_0$.
Suppose that $(\Omg, \mathcal{F}, \mu,\tau)$ is a Polish measure-preserving dynamical system equipped with its Borel $\sigma$-algebra and with a Borel probability measure. Suppose moreover that $\tau$ is surjective and aperiodic. 

\smallskip

Then there exists a non-trivial measurable partition $(A_1, A_2)$ of $\Omg$ with $\mu(A_1) > 1/2$ such that for the random products $T_n(\omg)$ given by
\[
T(\omega) =
\begin{cases}
	B_w & \text{if } \omega \in A_1,\\
	B_v & \text{if } \omega \in A_2,
\end{cases}
\] 
with 
\begin{align*}
	w_l = 1 + \frac{1}{l+1} \quad \textrm{and} \quad v_l =  1 - \frac{1}{l+1} \quad \textrm{for every $l \geq 1$}, 
\end{align*}
the sequence $(T_n(\omg))_{n \geq 1}$ is weakly mixing but not mixing for almost every $\omg \in \Omg$.
\epr

\smallskip

The proof of Proposition \ref{proppartitionharmonicsumsnotconverge} is based on \cite[Theorem 11]{DenielYves}, which is similar to \cite{akcogluJunco}, and relies on Rokhlin's lemma. This lemma was originally introduced for invertible aperiodic transformations. However, a version of this lemma also holds in the setting of a separable metric space endowed with its Borel $\sigma$-algebra and a Borel probability measure, together with a surjective and aperiodic transformation on this space.

\smallskip

More precisely, let $(\Omg,\mathcal{F},\mu,\tau)$ be a Polish measure-preserving dynamical system, where $\mathcal{F}$ is the Borel $\sigma$-algebra of $\Omg$ and $\mu$ is a Borel probability measure on this space. We suppose that $\tau$ is aperiodic, which means that for all $n \geq 1$, the set of points $x \in \Omg$ such that $\tau^n x = x$ has $\mu$-measure 0. We also suppose that $\tau$ is surjective. In this context, Rokhlin's lemma (see, for instance, \cite[Theorem 2.5 and Corollary 2.6]{Rokhlinlemmasurj}) states that for any $n \geq 1$ and for any $\varepsilon > 0$, there exists a measurable set $G \in \mathcal{F}$ such that the sets $G, \tau(G), \dotsc, \tau^{n-1}(G)$ are measurable and pairwise disjoint up to a set of measure 0 (which means that $\mu(\tau^i(G) \cap \tau^{j}(G)) = 0$ for every $0 \leq i \ne j \leq n-1$), $\mu(\tau^{l}(G)) = \mu(G)$ for every $0 \leq l \leq n-1$ and $\mu(\displaystyle \bigcup_{l=0}^{n-1} \tau^{l}(G)) > 1-\varepsilon$.

\smallskip

Note that on a non-atomic probability space, every ergodic transformation is aperiodic. Moreover, every Polish probability space equipped with its Borel $\sigma$-algebra and a continuous probability measure is isomorphic to $[0,1]$ equipped with the Lebesgue measure, see \cite[Theorem 17.41]{Ke}.

\smallskip

\bpf[Proof of Proposition \ref{proppartitionharmonicsumsnotconverge}]
Let $(n_j)_{j \geq 1}$ be a strictly increasing sequence of positive integers such that $\displaystyle \sum_{j \geq 1} 1/n_j^{1/3} < 1/3$. By Rokhlin's lemma, for every $j \geq 1$, there exists a measurable set $A_j$ such that $A_j, \dotsc, \tau^{n_j - 1} (A_j)$ are measurable, pairwise disjoint up to a set of measure 0, $\mu(\tau^l(A_j)) = \mu(A_j)$ for every $0 \leq l \leq n_j -1 $ and $\mu(\Omg \setminus \displaystyle \bigcup_{i=0}^{n_j - 1} \tau^i (A_j)) < 1/n_j^{1/3}$. Let us set $B_j = \displaystyle \bigcup_{i= n_j - n_j^{2/3}}^{n_j - 1} \tau^i (A_j)$ and $B = \displaystyle \bigcup_{j \geq 1} B_j$. Then 
$$
n_j^{1/3} \mu(B_j) = n_j^{1/3} n_j^{2/3} \mu(A_j) = \mu(\displaystyle \bigcup_{i=0}^{n_j - 1} \tau^i (A_j)) \leq 1
$$
and thus $\mu(B) \leq \displaystyle \sum_{j \geq 1} 1/n_j^{1/3} < 1/3$. We will now show that $\displaystyle \limsup_{n \to \infty} H_n \indic_B \geq 2/3$ for almost every $\omg \in \Omg$, where
$$
H_n \indic_B := \frac{1}{\log(n)} \displaystyle \sum_{i=0}^{n-1} \frac{1}{n-i} \indic_B \circ \tau^i.
$$
Let us set, for $j \geq 1$, $D_j = \displaystyle \bigcup_{i=0}^{n_j - n_j^{2/3} -1} \tau^i(A_j)$. Then
$$
\mu(D_j) + \mu(B_j) = \mu(\displaystyle \bigcup_{i=0}^{n_j - 1} \tau^{i}(A_j)) > 1 - 1/n_j^{1/3},
$$
and thus $\mu(D_j) > 1 - 2/n_j^{2/3}$.

Let $\omg \in D_j$ for $j \geq 1$, with $\omg \in \tau^i(A_j)$ and $0 \leq i \leq n_j - n_j^{2/3}-1$. Then $n_j - i \geq n_j^{2/3} + 1$ and
\begin{align*}
	H_{n_j - i} \indic_B(\omg) &= \frac{1}{\log(n_j - i)} \displaystyle \sum_{l=0}^{n_j - i -1} \frac{1}{n_j - i - l} \indic_B \circ \tau^l(\omg) \\
	&= \frac{1}{\log(n_j - i)} \displaystyle \sum_{l=1}^{n_j -i} \frac{1}{l} \indic_B \circ \tau^{n_j - i -l}(\omg)\\
	&\geq \frac{1}{\log(n_j - i)} \displaystyle \sum_{l=1}^{n_j^{2/3}} 1/l \\
	&\geq 2/3.
\end{align*}
Thus, for every $q \geq 1$, we have $\displaystyle \limsup_{n \to \infty} H_n \indic_B \geq 2/3$ on $\displaystyle \bigcap_{l \geq q} D_l$, and 
$$\mu \left(\displaystyle \bigcap_{l \geq q} D_l \right) \geq 1- \displaystyle \sum_{l\geq q} 2/n_l^{1/3}.$$ 
In particular, we have that $\displaystyle \limsup_{n \to \infty} H_n \indic_B(\omg) \geq 2/3$ for almost every $\omg \in \Omg$.

\smallskip

We now take $A_1 = \Omg \setminus B$, which satisfies $\mu(A_1) > 2/3$, and for almost every $\omg \in \Omg$, there exists a subsequence $(n_j)_{j \geq 1}$ such that 
$$
H_{n_j} \indic_{B}(\omg) \underset{j \to \infty}{\longrightarrow} \alpha
$$
with $\alpha \geq 2/3$.
In particular, we obtain that $H_{n_j} \indic_{A_1}(\omg)$ converges to $1-\alpha < 1/2$. The result follows from the fact that for almost every $\omg \in \Omg$, $$V_{n_j}(\omg) := \log \left(\displaystyle \prod_{l=1}^{n_j} \lvert \varepsilon_l(\tau^{n_j-l} \omg) \rvert \right) \underset{j \to \infty}{\longrightarrow} -\infty $$ 
by (\ref{ineqsupVn}).
\epf

\medskip

In the case where $\Omg = \T$, where $A_1,A_2$ are  two intervals of $[0,1)$, and where $\tau$ is an irrational rotation on $\Omg$, one can have more information on the dynamics of $(T_n(\omg))_{n \geq 1}$ for every $\omg \in \Omg$. It is based on an adaptation of Oxtoby's theorem for Nörlund sums. 

\smallskip

In what follows, we write $\mathcal{C}(\Omg)$ for the space of real continuous functions on $\Omg$ and $\lVert . \rVert_\infty$ for the supremum norm on $\mathcal{C}(\Omg)$.

\smallskip

\bpr \label{propOxtobynorlundmeans}
Let $(\Omg,\mathcal{F}, \mu)$ be a compact metric space equipped with its Borel $\sigma$-algebra and let $\mu$ be a Borel probability measure on $\Omg$. Let $\tau : \Omg \to \Omg$ be a continuous and uniquely ergodic transformation on $\Omg$ with respect to $\mu$. Let $(p_k)_{k \geq 1}$ be a sequence of positive real numbers such that 
\begin{align}
	&\frac{1}{P_n} (\lvert p_n - p_{n+1} \rvert + \dotsc + \lvert p_1 - p_{2} \rvert) \underset{n \to \infty}{\longrightarrow} 0, \\
	&\frac{p_{n+1}}{P_n} \underset{n \to \infty}{\longrightarrow} 0 \quad \textrm{and} \\
	&P_n \underset{n \to \infty}{\longrightarrow} \infty,
\end{align}
where $P_n = \displaystyle \sum_{j=1}^n p_j$.
Then for every continuous function $f : \Omg \to \R$, we have
\begin{align*}
	\frac{1}{P_n} \displaystyle \sum_{l=0}^{n-1} p_{n-l} f(\tau^l \omg) \underset{n \to \infty}{\longrightarrow} \int_\Omg f d\mu
\end{align*}
uniformly in $\omg \in \Omg$.
\epr 
\smallskip

\bpf
We adapt the proof given in \cite[Proposition 4.7.1]{BrinStuck}.
Suppose that the convergence does not hold. Then, there exists a continuous function $f : \Omg \to \R$, $\varepsilon > 0$, a strictly increasing sequence of positive integers $(n_j)_{j \geq 1}$ and a sequence $(x_j)_{j \geq 1} \subset \Omg$ such that for every $j \geq 1$, we have
\begin{align} \label{equationoxtobypoids}
	\left \lvert \frac{1}{P_{n_j}} \displaystyle \sum_{l=0}^{n_j-1} p_{n_j-l} f(\tau^l x_j) - \int_\Omg f d\mu  \right \rvert > \varepsilon.
\end{align}
Let us set
$$
U^{f}_{n_j}(x_j) :=  \frac{1}{P_{n_j}} \displaystyle \sum_{l=0}^{n_j-1} p_{n_j-l} f(\tau^l x_j).
$$
As in the Proof of \cite[Proposition 4.6.1]{BrinStuck}, one can show that there exists a subsequence $(n_{j_i})_{i \geq 1}$ of $(n_j)_{j \geq 1}$ such that $U_\infty(g) := \displaystyle \lim_{i \to \infty} U^{g}_{n_{j_i}}(x_{j_i}) $ exists for every $g \in C(\Omg)$ and defines a positive and bounded linear functional $L$ on $C(\Omg)$.
By Riesz's representation theorem, there exists a Borel probability measure $\nu$ on $\Omg$ such that
$$
U_\infty(g) = \int_\Omg g d\nu
$$    
for every $g \in C(\Omg)$. Let us notice that $U_\infty(g \circ \tau) = U_\infty(g)$. Indeed, for every $g \in C(\Omg)$, we have
\begin{align*}
	&\lvert U^{g \circ \tau}_{n_{j_i}}(x_{j_i}) - U^{g}_{n_{j_i}}(x_{j_i}) \rvert \\
	&= \left \lvert  \frac{1}{P_{n_{j_i}}} \displaystyle \sum_{l=0}^{n_{j_i}-1} p_{n_{j_i}-l} (g(\tau^{l+1} x_{j_i}) - g(\tau^{l} x_{j_i}) )   \right\rvert \\
	&=  \left \lvert  \frac{1}{P_{n_{j_i}}} \left(  \displaystyle \sum_{l=0}^{n_{j_i}-1} (p_{n_{j_i}-l +1} - p_{n_{j_i} -l}) g(\tau^{l} x_{j_i})     -p_{n_{j_i} +1} g(x_{j_i}) + p_1 g(\tau^{n_{j_i}} x_{j_i}) \right) \right\rvert \\
	&\leq \frac{1}{P_{n_{j_i}}}  \displaystyle \sum_{l=0}^{n_{j_i}-1} \lvert p_{n_{j_i}-l +1} - p_{n_{j_i} -l} \rvert \lVert g \rVert_\infty + \frac{1}{P_{n_{j_i}}} p_1 \lVert g \rVert_\infty + \frac{1}{P_{n_{j_i}}} p_{n_{j_i} +1} \lVert g \rVert_\infty
\end{align*}
which converges to 0 as $i \to \infty$ by assumption. Thus $U_\infty(g \circ \tau) = U_\infty(g)$ for every function $g \in C(\Omg)$, which means that $\nu$ is $\tau$-invariant. But this is a contradiction, since $\tau$ is uniquely ergodic for $\mu$ and
$$
\left \lvert \int_\Omg f d\nu - \int_\Omg f d\mu \right \rvert \geq \varepsilon
$$
by (\ref{equationoxtobypoids}).
\epf

\smallskip

Proposition \ref{propOxtobynorlundmeans} is also proved in \cite[Theorem 2]{Kozlov} in the case of an invertible map $\tau : \Omg \to \Omg$ and relies on \cite[Lemma 1 and Theorem 2, page 39]{SinaiCornfledFomin}, which asserts that the set $\{ g \circ \tau - g : g \in \mathcal{C}(\Omg)\}$ is dense in the set of functions $f \in \mathcal{C}(\Omg)$ satisfying $\int_\Omg f d\mu = 0$. This last fact uses that the transformation $\tau$ is invertible on $\Omg$, whereas our proof does not require it.

\smallskip
The result we obtain for two intervals of $[0,1)$ when $\Omg = \T$ is equipped with the normalized Lebesgue measure $m$ and an irrational rotation is the following.

\smallskip

\bpr
Let $\Omg = \T$ and let $m$ be the normalized Lebesgue measure on $\Omg$. Let $\tau  = R_\alpha$ with $\alpha \notin \Q$. Let $A_1, A_2$ be two disjoint intervals of $[0,1)$ such that $A_1 \cup A_2 = [0,1)$ and $\mu(A_l) > 0$ for $l = 1,2$.
Let 
\[
T(\omega) =
\begin{cases}
	B_w & \text{if } \omega \in A_1,\\
	B_v & \text{if } \omega \in A_2,
\end{cases}
\] 
where the weights are given by 
\begin{align*}
	w_l = 1 + \frac{1}{l+1} \quad \textrm{and} \quad v_l =  1 - \frac{1}{l+1} \quad \textrm{for every $l \geq 1$}. 
\end{align*}
\begin{enumerate}
	\item If $\mu(A_1) > \mu(A_2)$, then $(T_n(\omg))_{n \geq 1}$ is mixing for almost every $\omg \in \Omg$;
	\item If $\mu(A_1) < \mu(A_2)$, then $(T_n(\omg))_{n \geq 1}$ is not universal for almost every $\omg \in \Omg$.
\end{enumerate}
\epr

\smallskip

\bpf
Let us consider $f = \indic_{A_1}$. The weights are given by $p_n = 1/n $ for every $n \geq 1$, and we have that $P_n = \displaystyle \sum_{l=1}^{n} 1/l$ for every $n \geq 1$. 

For every $\varepsilon > 0$, there exist two real-valued continuous functions $f_1$ and $f_2$ on $\T$ such that $f_1 \leq f \leq f_2$ and $\int_{[0,1)} (f_2 - f_1) dm < \varepsilon$. From this, one can show using Proposition \ref{propOxtobynorlundmeans} that  \[\frac{1}{P_n} \displaystyle \sum_{l=0}^{n-1} \frac{1}{n-l} f(\tau^l x)  \underset{n \to \infty}{\longrightarrow} \mu(A_1)  \] 
for every $x \in \Omg$, and the proof follows exactly as in Proposition \ref{propmélangeexcasdeuxshiftsmesureA_1>1/2et<1/2}. Since the convergence holds for the whole sequence and not just for a subsequence, we conclude to the mixing property.
\epf

\medskip

It remains unknown of what happens in the case where $\mu(A_1) = \mu(A_2) = 1/2$, even if $\Omg = \T$, $\tau$ is the doubling map on $\Omg$, $A_1 = [0,1/2)$ and $A_2 = [1/2,1)$. In the latter case, the random variables $\indic_{A_1} \circ \tau^l$, $l \geq 0$, are independent and have the same distribution, namely $\nu := \frac{1}{2}(\delta_0 + \delta_1) $. Indeed, let $\varphi : \{0,1\}^{\Z_+} \to \T$ be given by $\varphi((x_j)_{j \geq 0}) = \displaystyle \sum_{j \geq 0} \frac{x_j}{2^{j+1}}$ and let $\sigma : (x_j)_{j \geq 0} \mapsto (x_{j+1})_{j \geq 0}$ be the backward shift on $\{0,1\}^{\Z_+}$. The map $\varphi$ is a bijection from the set $\Omg_0$ of binary sequences not eventually equal to $1$ onto $\T$, and $\tau \circ \varphi = \varphi \circ \sigma$. Let $\hat{X}_n$ be the projection of $\{0,1 \}^{\Z_+}$ onto the $n$-th coordinate. Let us notice that for every $l \geq 0$ and for every $\omg = \varphi(u) \in \T$, with $u \in \{0,1\}^{\Z_+}$, we have that $\indic_{A_1} \circ \tau^l(\omg) = \indic_{A_1}\circ \varphi \circ \sigma^l(u)$, which is equal to $1$ if and only if $\hat{X}_l(u) = 0$. Since the projections $\hat{X}_n, n \geq 0$, are independent on $\{0,1\}^{\Z_+}$ equipped with the $\sigma$-algebra generated by the cylinders and with the product measure $\nu^{\otimes \Z_+}$, it follows that the random variables $\indic_{A_1} \circ \tau^l$, $l \geq 0$, are independent and have the same distribution, which is $\nu$.

\smallskip
It would be natural to make use of central limit theorems in the case where $\mu(A_1) = 1/2$. However, the notion of a central limit theorem for Nörlund sums does not seem to have been developed in the setting of dynamical systems. By contrast, for independent random variables, one can refer to \cite{PeligradUtev}, which remains the most well-known work on the subject. Nevertheless, this result does not apply to our example in the case where $\Omega = \mathbb{T}$, $\tau$ is the doubling map, $A_1 = [0,1/2)$ and $A_2 = [1/2,1)$, since in this case, the variance of $\displaystyle \sum_{l=0}^{n-1} \frac{1}{n-l} (\indic_{[0,1/2)} -1/2)(\tau^l \omg)$ is equal to $\displaystyle \sum_{l=1}^{n} \frac{1}{l^2}$ and is bounded for $n \geq 1$.

\smallskip

Finally, we mention the paper by Chow and Lai \cite{ChowLai}, which investigates weighted sums similar to those considered in our example. While these results are technically applicable in our example when $\Omega = \mathbb{T}$, $\tau$ is the doubling map, $A_1 = [0,1/2)$ and $A_2 = [1/2,1)$, they do not provide any useful insight for our purposes.

\subsection{An example on the space $X = H(\C)$}
We consider in this subsection the case where $X = H(\C)$. We consider an ergodic measure-preserving dynamical system $(\Omg, \mathcal{F}, \mu, \tau)$, as well as a non-trivial measurable partition $(A_1,A_2)$ of $\Omg$. 
Let 
\[
T(\omega) =
\begin{cases}
	B_w & \text{if } \omega \in A_1,\\
	B_v & \text{if } \omega \in A_2,
\end{cases}
\] 
where the weights are given by 
\begin{align*}
	w_l = l \quad \textrm{and} \quad v_l =  1/l \quad \textrm{for every $l \geq 1$}. 
\end{align*}

In this case, we have to study
\begin{align}
	\log \left( \left(\displaystyle \prod_{l=1}^{n} \lvert \varepsilon_l(\tau^{n-l}\omg) \rvert \right) ^{1/n} \right) = \frac{1}{n} \displaystyle \sum_{l=1}^{n} \log(l) (\indic_{A_1}(\tau^{n-l} \omg) - \indic_{A_2}(\tau^{n-l} \omg)).
\end{align}
Let us observe that the sequence $(\log(n))_{n \geq 1}$ satisfies the assumptions of Lemma \ref{lemmecvgceL2sommesnorlund}. However, in the case of a strictly increasing sequence $(p_n)_{n \geq 1}$ of positive numbers such that $p_n/P_n \underset{n\to \infty}{\longrightarrow} 0$, where $P_n = \displaystyle \sum_{j=1}^n p_j$, it is known (see \cite[Theorem 20]{Hardybook}) that the Cesàro method is included in the Nörlund method for the positive weights $(p_n)_{n \geq 1}$. It means that if \(\displaystyle  \frac{1}{n} \sum_{l=1}^{n} s_l\) converges to $s \geq 0$, then \(\displaystyle \frac{1}{P_n} \sum_{l=1}^{n} p_l s_{n-l}\) also converges to $s$. Since $\displaystyle \sum_{l=1}^{n} \log(l) \sim n \log(n)$ as $n \to \infty$, an application of the Birkhoff ergodic theorem applied to the function $f := \indic_{A_1} - \indic_{A_2}$ leads to the following result.

\smallskip

\bpr \label{examplemelangeshiftslet1/lmesureA1diffmesureA2}
Let $X = H(\C)$.
Suppose that $(\Omg, \mathcal{F}, \mu,\tau)$ is an ergodic measure-preserving dynamical system and that $(A_1,A_2)$ is a partition of two measurable subsets of $\Omega$ with $\min \{\mu(A_1), \mu(A_2) \} > 0$.

Let 
\[
T(\omega) =
\begin{cases}
	B_w & \text{if } \omega \in A_1,\\
	B_v & \text{if } \omega \in A_2,
\end{cases}
\] 
where the weights are given by 
\begin{align*}
	w_l = l \quad \textrm{and} \quad v_l =  1/l \quad \textrm{for every $l \geq 1$}. 
\end{align*}
\begin{enumerate}
	\item If $\mu(A_1) > \mu(A_2)$, then $(T_n(\omg))_{n \geq 1}$ is mixing for almost every $\omg \in \Omg$;
	\item If $\mu(A_1) < \mu(A_2)$, then $(T_n(\omg))_{n \geq 1}$ is not universal for almost every $\omg \in \Omg$.
\end{enumerate}
\epr

\smallskip

In remains to study the case where $\mu(A_1) = \mu(A_2) = 1/2$. As in the previous example where $X = \ell_p$ or $X = c_0$, this case seems difficult to handle. However, the Central Limit Theorem from \cite{PeligradUtev} is applicable in the case of the doubling map and of $A_1 = [0,1/2)$ and $A_2 = [1/2,1)$. In fact, we can use a better result for our study, which is the following, due to Gut.

\smallskip

\bth[{\cite[Theorem 7.1]{gutarticleconvergence}}] \label{tharticleGutsommespoids}
Let $(Y_l)_{l \geq 0 }$ be a sequence of integrable i.i.d random variables such that $\mathbb{E}(Y_0) = 0$. 
Let $\{ (a_{n,l}, 0 \leq l \leq n-1), n \geq 1 \}$, be an array of positive real numbers.

Suppose that there exists $0 < p < 2$ such that $\mathbb{E}(\lvert Y_0 \rvert^{2p}) < \infty $. Suppose also that
\[
\displaystyle \sup_{\substack{n \geq 1 \\ 0 \leq l \leq n-1}} a_{n,l} < \infty.
\]
Then $n^{-1/p} \displaystyle \sum_{l=0}^{n-1} a_{n,l} Y_l$ converges to 0 almost surely.
\eth

\smallskip

Using Theorem \ref{tharticleGutsommespoids}, we can treat the situation where $\Omg = \T$, $\tau$ is the doubling map, $A_1 = [0,1/2)$ and $A_2 = [1/2,1)$.

\smallskip

\bpr \label{propfoncentierescasmuA_1=1/2doublangle}
Let $X = H(\C)$.
Suppose that $\Omg = \T$, $\tau$ is the doubling map on $\T$, $A_1 = [0,1/2)$ and $A_2 = [1/2,1)$.
Let 
\[
T(\omega) =
\begin{cases}
	B_w & \text{if } \omega \in A_1,\\
	B_v & \text{if } \omega \in A_2,
\end{cases}
\] 
where the weights are given by 
\begin{align*}
	w_l = l \quad \textrm{and} \quad v_l =  1/l \quad \textrm{for every $l \geq 1$}. 
\end{align*}
Then the sequence $(T_n(\omg))_{n \geq 1}$ is not universal, for almost every $\omg \in \Omg$.
\epr

\smallskip

\bpf
In the context of Proposition \ref{propfoncentierescasmuA_1=1/2doublangle}, the random variables $(\indic_{A_1} - \indic_{A_2})\circ \tau^l, l \geq 0$, are independent on $\T$. Moreover, applying Theorem \ref{tharticleGutsommespoids} with $p = 3/2$ and with $a_{n,l} = \frac{\log(n-l)}{n^{1/3}}$ for $0 \leq l \leq n-1$, we obtain that $\displaystyle \frac{1}{n} \displaystyle \sum_{l=0}^{n-1} \log(n-l) (\indic_{A_1} - \indic_{A_2})\circ \tau^l $ converges to 0 almost surely. In particular, the sequence $\left(\displaystyle \prod_{l=1}^{n} \lvert \varepsilon_l(\tau^{n-l}\omg) \rvert \right) ^{1/n}$ converges to $1$ for almost every $\omg \in \Omg$, and the sequence $(T_n(\omg))_{n \geq 1}$ is not universal for almost every $\omg \in \Omg$.
\epf

\medskip

\textbf{Acknowledgments.} The author gratefully thanks Quentin Menet for constructive feedback on a previous version of the manuscript, which contributed to refining the results and their presentation. He also thanks the referee for helpful comments.

\printbibliography

@book {BM,
	AUTHOR = {Bayart, Fr\'ed\'eric and Matheron, \'Etienne},
	TITLE = {Dynamics of linear operators},
	SERIES = {Cambridge Tracts in Mathematics},
	VOLUME = {179},
	PUBLISHER = {Cambridge University Press, Cambridge},
	YEAR = {2009},
	PAGES = {xiv+337},
	ISBN = {978-0-521-51496-5},
	MRCLASS = {47-02 (11M06 37B05 47A16 47A35)},
	MRNUMBER = {2533318},
	MRREVIEWER = {E.\ A.\ Gallardo-Guti\'errez},
	DOI = {10.1017/CBO9780511581113},
	URL = {https://doi.org/10.1017/CBO9780511581113},
}

@article {E,
	AUTHOR = {Eisner, Tanja},
	TITLE = {A ``typical'' contraction is unitary},
	JOURNAL = {Enseign. Math. (2)},
	FJOURNAL = {L'Enseignement Math\'ematique. Revue Internationale. 2e
	S\'erie},
	VOLUME = {56},
	YEAR = {2010},
	NUMBER = {3-4},
	PAGES = {403--410},
	ISSN = {0013-8584},
	MRCLASS = {47A35 (46C05 47D06 47L99)},
	MRNUMBER = {2769030},
	MRREVIEWER = {Wojciech\ Bartoszek},
	DOI = {10.4171/LEM/56-3-6},
	URL = {https://doi.org/10.4171/LEM/56-3-6},
}

@misc{Gill3,
%	title={Linear dynamics of random products of operators}, 
%	author={Valentin Gillet},
%	year={2025},
%	eprint={2507.00186},
%	archivePrefix={arXiv},
%	primaryClass={math.FA},
%	url={https://arxiv.org/abs/2507.00186}, 
%}

@book {GEP,
	AUTHOR = {Grosse-Erdmann, Karl-G. and Peris Manguillot, Alfredo},
	TITLE = {Linear chaos},
	SERIES = {Universitext},
	PUBLISHER = {Springer, London},
	YEAR = {2011},
	PAGES = {xii+386},
	ISBN = {978-1-4471-2169-5},
	MRCLASS = {47-02 (37D45 47A16)},
	MRNUMBER = {2919812},
	MRREVIEWER = {E.\ A.\ Gallardo-Guti\'errez},
	DOI = {10.1007/978-1-4471-2170-1},
	URL = {https://doi.org/10.1007/978-1-4471-2170-1},
}

@book {Ke,
	AUTHOR = {Kechris, Alexander S.},
	TITLE = {Classical descriptive set theory},
	SERIES = {Graduate Texts in Mathematics},
	VOLUME = {156},
	PUBLISHER = {Springer-Verlag, New York},
	YEAR = {1995},
	PAGES = {xviii+402},
	ISBN = {0-387-94374-9},
	MRCLASS = {03E15 (03-01 03-02 04A15 28A05 54H05 90D44)},
	MRNUMBER = {1321597},
	MRREVIEWER = {Jakub\ Jasi\'nski},
	DOI = {10.1007/978-1-4612-4190-4},
	URL = {https://doi.org/10.1007/978-1-4612-4190-4},
}

@article {L,
	AUTHOR = {Lomonosov, V. I.},
	TITLE = {Invariant subspaces of the family of operators that commute
	with a completely continuous operator},
	JOURNAL = {Funkcional. Anal. i Prilo\v zen.},
	FJOURNAL = {Akademija Nauk SSSR. Funkcional\cprime nyi Analiz i ego
	Prilo\v zenija},
	VOLUME = {7},
	YEAR = {1973},
	NUMBER = {3},
	PAGES = {55--56},
	ISSN = {0374-1990},
	MRCLASS = {47A15},
	MRNUMBER = {420305},
	MRREVIEWER = {R.\ G.\ Bartle},
}

@article {Bell,
	AUTHOR = {Bellman, Richard},
	TITLE = {Limit theorems for non-commutative operations. {I}},
	JOURNAL = {Duke Math. J.},
	FJOURNAL = {Duke Mathematical Journal},
	VOLUME = {21},
	YEAR = {1954},
	PAGES = {491--500},
	ISSN = {0012-7094,1547-7398},
	MRCLASS = {60.0X},
	MRNUMBER = {62368},
	MRREVIEWER = {J.\ Wolfowitz},
	URL = {http://projecteuclid.org/euclid.dmj/1077465878},
}

@article {CB,
	AUTHOR = {Conze, Jean-Pierre and Le Borgne, St\'ephane},
	TITLE = {On the {CLT} for rotations and {BV} functions},
	JOURNAL = {Ann. Math. Blaise Pascal},
	FJOURNAL = {Annales Math\'ematiques Blaise Pascal},
	VOLUME = {29},
	YEAR = {2022},
	NUMBER = {1},
	PAGES = {51--97},
	ISSN = {1259-1734,2118-7436},
	MRCLASS = {11K50 (11K60 37E10 60F05)},
	MRNUMBER = {4490762},
	MRREVIEWER = {Thomas\ Ward},
	DOI = {10.5802/ambp.407},
	URL = {https://doi.org/10.5802/ambp.407},
}

@article {FK,
	AUTHOR = {Furstenberg, H. and Kesten, H.},
	TITLE = {Products of random matrices},
	JOURNAL = {Ann. Math. Statist.},
	FJOURNAL = {Annals of Mathematical Statistics},
	VOLUME = {31},
	YEAR = {1960},
	PAGES = {457--469},
	ISSN = {0003-4851},
	MRCLASS = {60.00},
	MRNUMBER = {121828},
	MRREVIEWER = {R.\ E.\ Bellman},
	DOI = {10.1214/aoms/1177705909},
	URL = {https://doi.org/10.1214/aoms/1177705909},
}

@article {Ha,
	AUTHOR = {Hal\'asz, G.},
	TITLE = {Remarks on the remainder in {B}irkhoff's ergodic theorem},
	JOURNAL = {Acta Math. Acad. Sci. Hungar.},
	FJOURNAL = {Acta Mathematica. Academiae Scientiarum Hungaricae},
	VOLUME = {28},
	YEAR = {1976},
	NUMBER = {3-4},
	PAGES = {389--395},
	ISSN = {0001-5954,1588-2632},
	MRCLASS = {28A65},
	MRNUMBER = {425076},
	MRREVIEWER = {Meir\ Smorodinsky},
	DOI = {10.1007/BF01896805},
	URL = {https://doi.org/10.1007/BF01896805},
}

@article {K,
	AUTHOR = {Kac, M.},
	TITLE = {On the distribution of values of sums of the type {$\sum f(2^k
	t)$}},
	JOURNAL = {Ann. of Math. (2)},
	FJOURNAL = {Annals of Mathematics. Second Series},
	VOLUME = {47},
	YEAR = {1946},
	PAGES = {33--49},
	ISSN = {0003-486X},
	MRCLASS = {42.4X},
	MRNUMBER = {15548},
	MRREVIEWER = {A.\ Zygmund},
	DOI = {10.2307/1969033},
	URL = {https://doi.org/10.2307/1969033},
}

@article {Kachu,
	AUTHOR = {Kachurovski\u i, A. G.},
	TITLE = {Rates of convergence in ergodic theorems},
	JOURNAL = {Uspekhi Mat. Nauk},
	FJOURNAL = {Uspekhi Matematicheskikh Nauk},
	VOLUME = {51},
	YEAR = {1996},
	NUMBER = {4(310)},
	PAGES = {73--124},
	ISSN = {0042-1316,2305-2872},
	MRCLASS = {28D05 (47A35 60F15 60G10)},
	MRNUMBER = {1422228},
	MRREVIEWER = {Manfred\ Denker},
	%DOI = {10.1070/RM1996v051n04ABEH002964},
	URL = {https://doi.org/10.1070/RM1996v051n04ABEH002964},
}

@book {N,
	AUTHOR = {Nikolski, Nikola\"i},
	TITLE = {Hardy spaces},
	SERIES = {Cambridge Studies in Advanced Mathematics},
	VOLUME = {179},
	EDITION = {French},
	PUBLISHER = {Cambridge University Press, Cambridge},
	YEAR = {2019},
	PAGES = {xviii+277},
	ISBN = {978-1-107-18454-1},
	MRCLASS = {42-01 (42B30 42B35)},
	MRNUMBER = {3890074},
}

@book {Pet,
	AUTHOR = {Petersen, Karl},
	TITLE = {Ergodic theory},
	SERIES = {Cambridge Studies in Advanced Mathematics},
	VOLUME = {2},
	PUBLISHER = {Cambridge University Press, Cambridge},
	YEAR = {1983},
	PAGES = {xii+329},
	ISBN = {0-521-23632-0},
	MRCLASS = {28-02 (28D05 54H20 58F11)},
	MRNUMBER = {833286},
	MRREVIEWER = {Nathaniel\ F. G. Martin},
	DOI = {10.1017/CBO9780511608728},
	URL = {https://doi.org/10.1017/CBO9780511608728},
}

@article {S,
	AUTHOR = {Schmidt, Klaus},
	TITLE = {On recurrence},
	JOURNAL = {Z. Wahrsch. Verw. Gebiete},
	FJOURNAL = {Zeitschrift f\"ur Wahrscheinlichkeitstheorie und Verwandte
	Gebiete},
	VOLUME = {68},
	YEAR = {1984},
	NUMBER = {1},
	PAGES = {75--95},
	ISSN = {0044-3719},
	MRCLASS = {60F15 (28D99 60G10)},
	MRNUMBER = {767446},
	MRREVIEWER = {Eric\ V.\ Slud},
	DOI = {10.1007/BF00535175},
	URL = {https://doi.org/10.1007/BF00535175},
}

@book {Wal,
	AUTHOR = {Walters, Peter},
	TITLE = {An introduction to ergodic theory},
	SERIES = {Graduate Texts in Mathematics},
	VOLUME = {79},
	PUBLISHER = {Springer-Verlag, New York-Berlin},
	YEAR = {1982},
	PAGES = {ix+250},
	ISBN = {0-387-90599-5},
	MRCLASS = {28Dxx (54H20 58F11)},
	MRNUMBER = {648108},
	MRREVIEWER = {M.\ A.\ Akcoglu},
}

@book {BrinStuck,
	AUTHOR = {Brin, Michael and Stuck, Garrett},
	TITLE = {Introduction to dynamical systems},
	PUBLISHER = {Cambridge University Press, Cambridge},
	YEAR = {2015},
	PAGES = {xii+247},
	ISBN = {978-1-107-53894-8; 978-0-521-80841-5},
	MRCLASS = {37-01},
	MRNUMBER = {3558919},
}

@article {Kozlov,
	AUTHOR = {Kozlov, V. V.},
	TITLE = {Weighted Means, Strict Ergodicity, and Uniform Distributions
	},
	JOURNAL = {Mathematical Notes},
	VOLUME = {78},
	YEAR = {2005},
	NUMBER = {3},
	PAGES = {329–337},
	DOI = {10.1007/s11006-005-0132-x},
	URL = {https://doi.org/10.1007/s11006-005-0132-x},
}

@book {SinaiCornfledFomin,
	AUTHOR = {Cornfeld, I. P. and Fomin, S. V. and Sinai, Ya. G.},
	TITLE = {Ergodic theory},
	SERIES = {Grundlehren der mathematischen Wissenschaften [Fundamental
	Principles of Mathematical Sciences]},
	VOLUME = {245},
	NOTE = {Translated from the Russian by A. B. Sosinskii},
	PUBLISHER = {Springer-Verlag, New York},
	YEAR = {1982},
	PAGES = {x+486},
	ISBN = {0-387-90580-4},
	MRCLASS = {28D05 (54H20 58F11)},
	MRNUMBER = {832433},
	DOI = {10.1007/978-1-4615-6927-5},
	URL = {https://doi.org/10.1007/978-1-4615-6927-5},
}

@article {DenielYves,
	AUTHOR = {D\'eniel, Yves},
	TITLE = {On the a.s.\ {C}es\`aro-{$\alpha$} convergence for stationary
	or orthogonal random variables},
	JOURNAL = {J. Theoret. Probab.},
	FJOURNAL = {Journal of Theoretical Probability},
	VOLUME = {2},
	YEAR = {1989},
	NUMBER = {4},
	PAGES = {475--485},
	ISSN = {0894-9840,1572-9230},
	MRCLASS = {28D05 (42C15)},
	MRNUMBER = {1011200},
	MRREVIEWER = {U.\ Krengel},
	DOI = {10.1007/BF01051879},
	URL = {https://doi.org/10.1007/BF01051879},
}

@article{ChowLai,
	ISSN = {00911798, 2168894X},
	URL = {http://www.jstor.org/stable/2959450},
	author = {Y. S. Chow and T. L. Lai},
	journal = {The Annals of Probability},
	number = {5},
	pages = {810--824},
	publisher = {Institute of Mathematical Statistics},
	title = {Limiting Behavior of Weighted Sums of Independent Random Variables},
	urldate = {2025-10-12},
	volume = {1},
	year = {1973}
}

@article{PeligradUtev,
	author = {Magda Peligrad and Sergey Utev},
	title = {{Central limit theorem for linear processes}},
	volume = {25},
	journal = {The Annals of Probability},
	number = {1},
	publisher = {Institute of Mathematical Statistics},
	pages = {443 -- 456},
	year = {1997},
	doi = {10.1214/aop/1024404295},
	URL = {https://doi.org/10.1214/aop/1024404295}
}

@book {Hardybook,
	AUTHOR = {Hardy, G. H.},
	TITLE = {Divergent series},
	NOTE = {With a preface by J. E. Littlewood and a note by L. S.
	Bosanquet,
	Reprint of the revised (1963) edition},
	PUBLISHER = {\'Editions Jacques Gabay, Sceaux},
	YEAR = {1992},
	PAGES = {xvi+396},
	ISBN = {2-87647-131-0},
	MRCLASS = {01A75 (40-01 40-03)},
	MRNUMBER = {1188874},
}

@article {gutarticleconvergence,
	AUTHOR = {Gut, A.},
	TITLE = {Complete convergence for arrays},
	JOURNAL = {Period. Math. Hungar.},
	FJOURNAL = {Periodica Mathematica Hungarica. Journal of the J\'anos Bolyai
	Mathematical Society},
	VOLUME = {25},
	YEAR = {1992},
	NUMBER = {1},
	PAGES = {51--75},
	ISSN = {0031-5303,1588-2829},
	MRCLASS = {60F15 (60G50)},
	MRNUMBER = {1200841},
	MRREVIEWER = {Zdzis\l aw\ Rychlik},
	DOI = {10.1007/BF02454383},
	URL = {https://doi.org/10.1007/BF02454383},
}

@article {akcogluJunco,
	AUTHOR = {Akcoglu, M. A. and del Junco, A.},
	TITLE = {Convergence of averages of point transformations},
	JOURNAL = {Proc. Amer. Math. Soc.},
	FJOURNAL = {Proceedings of the American Mathematical Society},
	VOLUME = {49},
	YEAR = {1975},
	PAGES = {265--266},
	ISSN = {0002-9939,1088-6826},
	MRCLASS = {28A65},
	MRNUMBER = {360999},
	MRREVIEWER = {Michael\ Lin},
	DOI = {10.2307/2039829},
	URL = {https://doi.org/10.2307/2039829},
}

@article {Rokhlinlemmasurj,
	AUTHOR = {Heinemann, Stefan-M. and Schmitt, Oliver},
	TITLE = {Rokhlin's lemma for non-invertible maps},
	JOURNAL = {Dynam. Systems Appl.},
	FJOURNAL = {Dynamic Systems and Applications},
	VOLUME = {10},
	YEAR = {2001},
	NUMBER = {2},
	PAGES = {201--213},
	ISSN = {1056-2176},
	MRCLASS = {28D05 (37A05 43A05)},
	MRNUMBER = {1843737},
	MRREVIEWER = {Idris Assani},
}

@article {Bayart-Costakis-Hadjiloucas10,
	AUTHOR = {Bayart, Fr\'ed\'eric and Costakis, George and Hadjiloucas,
	Demetris},
	TITLE = {Topologically transitive skew-products of operators},
	JOURNAL = {Ergodic Theory Dynam. Systems},
	FJOURNAL = {Ergodic Theory and Dynamical Systems},
	VOLUME = {30},
	YEAR = {2010},
	NUMBER = {1},
	PAGES = {33--49},
	ISSN = {0143-3857,1469-4417},
	MRCLASS = {47A16 (37B05 47B37)},
	MRNUMBER = {2586344},
	MRREVIEWER = {Enhui\ Shi},
	DOI = {10.1017/S0143385708001065},
	URL = {https://doi.org/10.1017/S0143385708001065},
}

@article {Costakis-Hadjiloucas08,
	AUTHOR = {Costakis, George and Hadjiloucas, Demetris},
	TITLE = {Topologically transitive skew-products of backward shift
	operators and hypercyclicity},
	JOURNAL = {Proc. Amer. Math. Soc.},
	FJOURNAL = {Proceedings of the American Mathematical Society},
	VOLUME = {136},
	YEAR = {2008},
	NUMBER = {3},
	PAGES = {937--946},
	ISSN = {0002-9939,1088-6826},
	MRCLASS = {47A16 (28D05 37B05 47A35)},
	MRNUMBER = {2361867},
	MRREVIEWER = {Sophie\ Grivaux},
	DOI = {10.1090/S0002-9939-07-09001-6},
	URL = {https://doi.org/10.1090/S0002-9939-07-09001-6},
}

@article {Martinez-Mendez-Peris26,
	AUTHOR = {Mart\'inez-Gim\'enez, F\'elix and M\'endez-G\'omez, H\'ector
	and Peris, Alfred and Rodenas, Francisco},
	TITLE = {Chaotic skew-products of operators on {F}r\'echet spaces},
	JOURNAL = {J. Math. Anal. Appl.},
	FJOURNAL = {Journal of Mathematical Analysis and Applications},
	VOLUME = {553},
	YEAR = {2026},
	NUMBER = {2},
	PAGES = {Paper No. 129949, 11},
	ISSN = {0022-247X,1096-0813},
	MRCLASS = {47A16 (30H10 37A25)},
	MRNUMBER = {4944038},
	MRREVIEWER = {Epaminondas\ A.\ Diamantopoulos},
	DOI = {10.1016/j.jmaa.2025.129949},
	URL = {https://doi.org/10.1016/j.jmaa.2025.129949},
}

@article {Godefroy-Shapiro91,
	AUTHOR = {Godefroy, Gilles and Shapiro, Joel H.},
	TITLE = {Operators with dense, invariant, cyclic vector manifolds},
	JOURNAL = {J. Funct. Anal.},
	FJOURNAL = {Journal of Functional Analysis},
	VOLUME = {98},
	YEAR = {1991},
	NUMBER = {2},
	PAGES = {229--269},
	ISSN = {0022-1236,1096-0783},
	MRCLASS = {47A65 (47B37 47B38 58F13)},
	MRNUMBER = {1111569},
	MRREVIEWER = {Warren\ R.\ Wogen},
	DOI = {10.1016/0022-1236(91)90078-J},
	URL = {https://doi.org/10.1016/0022-1236(91)90078-J},
}

@article {Gonzalez-Quas14,
	AUTHOR = {Gonz\'alez-Tokman, Cecilia and Quas, Anthony},
	TITLE = {A semi-invertible operator {O}seledets theorem},
	JOURNAL = {Ergodic Theory Dynam. Systems},
	FJOURNAL = {Ergodic Theory and Dynamical Systems},
	VOLUME = {34},
	YEAR = {2014},
	NUMBER = {4},
	PAGES = {1230--1272},
	ISSN = {0143-3857,1469-4417},
	MRCLASS = {47A35 (37A30 37H15)},
	MRNUMBER = {3227155},
	MRREVIEWER = {Wojciech\ Bartoszek},
	DOI = {10.1017/etds.2012.189},
	URL = {https://doi.org/10.1017/etds.2012.189},
}

@article {Froyland-Lloyd-Quas13,
	AUTHOR = {Froyland, Gary and Lloyd, Simon and Quas, Anthony},
	TITLE = {A semi-invertible {O}seledets theorem with applications to
	transfer operator cocycles},
	JOURNAL = {Discrete Contin. Dyn. Syst.},
	FJOURNAL = {Discrete and Continuous Dynamical Systems},
	VOLUME = {33},
	YEAR = {2013},
	NUMBER = {9},
	PAGES = {3835--3860},
	ISSN = {1078-0947,1553-5231},
	MRCLASS = {37H15 (37A30 37L55)},
	MRNUMBER = {3038042},
	MRREVIEWER = {Luciana\ A.\ Alves},
	DOI = {10.3934/dcds.2013.33.3835},
	URL = {https://doi.org/10.3934/dcds.2013.33.3835},
}

\end{document}